\documentclass[oneside, reqno]{amsart}

\usepackage{xypic}
\input xy
\xyoption{all}
\usepackage{epsfig}
\usepackage{amsthm}
\usepackage{amssymb}
\usepackage{amsmath}
\usepackage{amscd}
\usepackage{color}
\usepackage[T1]{fontenc}
\usepackage{stackengine}
\usepackage{upgreek}
\usepackage[font=scriptsize]{caption}
\usepackage{wrapfig}
\stackMath

\newcommand{\bg}{\begin{equation}}
\newcommand{\ed}{\end{equation}}
\newcommand{\bga}{\begin{eqnarray}}
\newcommand{\eda}{\end{eqnarray}}
\newcommand{\pf}{\textbf{Proof:\ }}

\def\cbdu{\par{\raggedleft$\Box$\par}}

\newtheorem {Theorem}  {Theorem}

\numberwithin{Theorem}{section}

\newtheorem {Lemma}[Theorem]  {Lemma}
\newtheorem {Proposition}[Theorem]{Proposition}
\theoremstyle{definition}
\newtheorem{Definition}[Theorem]{Definition}
\theoremstyle{remark}

\expandafter\chardef\csname pre amssym.def
at\endcsname=\the\catcode`\@ \catcode`\@=11
\def\undefine#1{\let#1\undefined}
\def\newsymbol#1#2#3#4#5{\let\next@\relax
 \ifnum#2=\@ne\let\next@\msafam@\else
 \ifnum#2=\tw@\let\next@\msbfam@\fi\fi
 \mathchardef#1="#3\next@#4#5}
\def\mathhexbox@#1#2#3{\relax
 \ifmmode\mathpalette{}{\m@th\mathchar"#1#2#3}%
 \else\leavevmode\hbox{$\m@th\mathchar"#1#2#3$}\fi}
\def\hexnumber@#1{\ifcase#1 0\or 1\or 2\or 3\or 4\or 5\or 6\or 7\or 8\or
 9\or A\or B\or C\or D\or E\or F\fi}

\font\teneufm=eufm10 \font\seveneufm=eufm7 \font\fiveeufm=eufm5
\newfam\eufmfam
\textfont\eufmfam=\teneufm \scriptfont\eufmfam=\seveneufm
\scriptscriptfont\eufmfam=\fiveeufm

\catcode`\@=\csname pre amssym.def at\endcsname

\newcounter{remark}
\def  \12  {{\frac{1}{2}}}

\def\build#1_#2^#3{\mathrel{\mathop{\kern 0pt#1}\limits_{#2}^{#3}}}

\numberwithin{equation}{section}

\begin{document}

\title[Almost sure existence for EMHD]{Almost sure existence of global weak solutions for supercritical electron MHD}


\author [Mimi Dai]{Mimi Dai}

\address{Department of Mathematics, Statistics and Computer Science, University of Illinois at Chicago, Chicago, IL 60607, USA}
\address{School of Mathematics, Institute for Advanced Study, Princeton, NJ 08540, USA}
\email{mdai@uic.edu}


\begin{abstract}
We consider the Cauchy problem for the electron magnetohydrodynamics model in the supercritical regime. For rough initial data in $\mathcal H^{-s}(\mathbb T^n)$ with $s>0$, we obtain global in time weak solutions almost surely via an appropriate randomization of the initial data.

\bigskip


\vspace{0.05cm}
CLASSIFICATION CODE: 35Q35, 76D03, 76W05.
\end{abstract}

\maketitle

\section{Introduction}

With static background flow the electron magnetohydrodynamics(MHD) is modeled by 
\begin{equation}\label{emhd0}
\begin{split}
B_t+\nabla\times((\nabla\times B)\times B)=&\ \Delta B,\\
\nabla\cdot B=&\ 0
\end{split}
\end{equation}
where the nonlinear term captures the Hall effect. It is a subsystem of the full MHD system with Hall effect, which has attracted much attention in recent time, see \cite{ADFL, CDL, CWeng}. The magnetic field considered in this paper takes the form $B(x,t)=(B_1(x, t), B_2(x, t),B_3(x, t))$ for either $x\in \mathbb T^2$ or $x\in \mathbb T^3$. In particular, in the former case, the problem is regarded as for a 3-dimensional (3D) magnetic field posed on the 2D torus (c.f. \cite{CW}).
One notices that the highest order derivative appears both in the linear diffusion and the quadratic Hall term, resulting (\ref{emhd0}) as a quasilinear system. Beside other difficulties caused by the peculiar geometry structure of the Hall effect, the quasilinear feature is a major obstacle in the analysis of (\ref{emhd0}). We will further illustrate this point by the discussion of scaling property. System (\ref{emhd0}) has the natural scaling that if $B(x,t)$ solves the system with initial data $B_0(x)$, then the rescaled vector field
\[B_\lambda(x,t)=B(\lambda x, \lambda^{2}t)\]
for an arbitrary parameter $\lambda$
solves the system as well with initial data $B_0(\lambda x)$. Some scaling invariant spaces (also referred as critical spaces) with embedding for (\ref{emhd0}) in $n$-dimensional space are
\begin{equation}\label{critical1}
\dot{\mathcal H}^{\frac{n}2}\hookrightarrow L^{\infty}\hookrightarrow \dot B^{\frac{n}{p}}_{p,\infty} \hookrightarrow BMO \hookrightarrow \dot B^{0}_{\infty,\infty}, \ \ \ 1<p<\infty.
\end{equation} 
In view of (\ref{critical1}) we note the energy space $L^2(\mathbb T^n)$ is supercritical in both 2D and 3D, and system (\ref{emhd0}) is supercritical in both situations. Thus it is naturally challenging to analyze (\ref{emhd0}) even in 2D.

From the perspective of mathematics, in order to understand the competition between the nonlinear Hall effect and the linear term, we consider the electron MHD with generalized diffusion
\begin{equation}\label{emhd}
\begin{split}
B_t+\nabla\times((\nabla\times B)\times B)=&-(-\Delta)^\alpha B,\\
\nabla\cdot B=&\ 0
\end{split}
\end{equation}
on $[0,\infty)\times \mathbb T^n$, $n=2,3$, for $\alpha>0$. The MHD system with fractional diffusion was previously studied, for instance see \cite{CWW}. System (\ref{emhd}) possesses the scaling
\begin{equation}\notag
B_\lambda(x,t)=\lambda^{2\alpha-2}B(\lambda x, \lambda^{2\alpha}t),
\end{equation}
according to which, some critical spaces with embedding for (\ref{emhd}) are
\begin{equation}\label{critical2}
\dot{\mathcal H}^{\frac{n}2+2-2\alpha}\hookrightarrow L^{\frac{n}{2\alpha-2}}\hookrightarrow \dot B^{2-2\alpha+\frac{n}{p}}_{p,\infty} \hookrightarrow BMO^{2-2\alpha} \hookrightarrow \dot B^{2-2\alpha}_{\infty,\infty}, \ 1<p<\infty.
\end{equation} 
The basic energy law for (\ref{emhd}) is given by
\begin{equation}\label{energy}
\|B(x,t)\|_{L^2}^2+\int_0^t\|\nabla^\alpha B(x,s)\|_{L^2}^2\, ds=\|B(x,0)\|_{L^2}^2.
\end{equation}
It follows from (\ref{energy}) that solutions of (\ref{emhd}) satisfy the a priori estimates 
\[B\in L^\infty([0,T); L^2(\mathbb T^n)) \cap L^2([0,T); H^\alpha (\mathbb T^n)).\]
On the other hand we see from (\ref{critical2}) that the Sobolev space $\dot{\mathcal H}^{3-2\alpha}$ is critical for (\ref{emhd}) in 2D, while $\dot{\mathcal H}^{\frac72-2\alpha}$ is critical in 3D.  Since $\dot{\mathcal H}^{3-2\alpha}=L^2$ for $\alpha=\frac32$ and $\dot{\mathcal H}^{\frac72-2\alpha}=L^2$ for $\alpha=\frac74$, system (\ref{emhd}) is critical in 2D when $\alpha=\frac32$ and critical in 3D when $\alpha=\frac74$.  For $\alpha>\frac32$ in 2D and $\alpha>\frac74$ in 3D, system (\ref{emhd}) is referred to be subcritical in which situation the linear term dominates, and hence global regularity is known to hold by standard energy method. While for $\alpha<\frac32$ in 2D and $\alpha<\frac74$ in 3D (including $\alpha=1$), system (\ref{emhd}) is supercritical and challenging in general. In this paper, we consider (\ref{emhd}) in the supercritical and critical setting, i.e. $\alpha\leq\frac32$ in 2D and $\alpha\leq\frac74$ in 3D.

 When the initial data is rather regular, say in $\dot{\mathcal H}^{s}$ with $s>\frac n2+2-2\alpha$, well-posedness of (\ref{emhd}) is expected, see \cite{Dai2} for instance. With initial data $B_0\in L^2(\mathbb T^n)$, one can obtain weak solutions for (\ref{emhd}) by using Galerkin approximating approach.
Nevertheless, for rough initial data below $L^2(\mathbb T^n)$, it is not clear how to construct weak solutions for (\ref{emhd}). This is a similar situation for many other equations, like the Navier-Stokes equation, nonlinear Schr\"ondinger equation, nonlinear wave equation, etc. 
The purpose of the paper is to construct global in time weak solutions to (\ref{emhd}) by randomizing initial data in Sobolev spaces $\dot{\mathcal H}^{-s}(\mathbb T^n)$ with $n=2,3$ and $s>0$. For a given rough initial data $f\in \dot{\mathcal H}^{-s}(\mathbb T^n)$ with $\nabla\cdot f=0$ and $\int_{\mathbb T^n} f\, dx=0$, we randomize it appropriately to $f^\omega$ satisfying $\nabla\cdot f^\omega=0$. We then consider solution of (\ref{emhd}) with the initial data $f^\omega$ in the form
\[B(x,t)= e^{-t(-\Delta)^\alpha}f^\omega (x)+H(x,t)\]
where $H$ satisfies a nonlinear equation that depends on $e^{-t(-\Delta)^\alpha}f^\omega$ and obviously $H(x,0)=0$. A crucial point is that the free evolution $e^{-t(-\Delta)^\alpha}f^\omega$ has almost surely improved $L^p$ estimates thanks to the randomization of the data. As a consequence, it provides the possibility to construct a global in time weak solution $H$. 

The study of well-posedness for randomized initial data was initiated by Bourgain in \cite{Bou96} for supercritical nonlinear Schr\"ondinger equation. Random data Cauchy problem was investigated for supercritical wave equation by Burq and Tzvetkov \cite{BT1, BT2}, eventually leading to a global existence theory. Applying the method from \cite{BT1, BT2}, Nahmod, Pavlovi\'c and Staffilani \cite{NPS} showed almost sure existence of global weak solutions for the Navier-Stokes equation (NSE) below $L^2(\mathbb T^n)$. An almost sure global well-posedness result for the 2D nonlinear Schr\"ondinger equation with random radial initial data in the supercritical regime was established by Deng \cite{Deng}.
Later on, L\"uhrmann and Mendelson \cite{LM} established the random data Cauchy theory for nonlinear wave equations of power-type on $\mathbb R^3$. 
Random data Cauchy problem has been also treated for various other equations when the deterministic Cauchy theory is hard to be achieved. We do not intend to give an extensive list here. 

Applying the framework of random data Cauchy theory to (\ref{emhd}) in this paper, the main difficulty comes from the strong nonlinear effect. The Hall term of (\ref{emhd}) is one degree higher than the nonlinear term $(u\cdot\nabla) u$ of the well-known (NSE). It is a general belief that the nonlinearity of $(u\cdot\nabla) u$ poses intrinsic obstacles to crack the global regularity problem for the 3D NSE, see \cite{Tao}. 
As a result of the presence of the strong nonlinear term in (\ref{emhd}) and the quasilinear feature of (\ref{emhd}), it prevents us to show global existence of weak solutions with randomized data for $\alpha=1$ in both 2D and 3D. Nevertheless, for larger value of $\alpha$ which is still below the critical exponent, we are able to obtain almost sure existence of global weak solutions for (\ref{emhd}) in the supercritical regime.

\bigskip

\section{Main results}\label{statement}

In this section we fix notations to be used throughout the text, introduce the procedure of randomization, and then state the main results. 

\medskip

\subsection{Notations}
We often denote $C$ by a constant in estimates which may vary from line to line.
When it is not necessary to track the constant, $f\lesssim g$ is used to denote $f\leq Cg$ for some constant $C>0$.

Note the inhomogeneous and homogeneous Sobolev spaces are equivalent on torus. In the rest of the paper we only use $\mathcal H^s$ to denote the Sobolev space. We further denote
\begin{equation}\notag
\begin{split}
\mathcal H=&\ \mbox{the closure of} \ \{ f\in C^\infty(\mathbb T^n)| \nabla\cdot f=0\} \ \mbox{in} \ L^2(\mathbb T^n),\\
\mathcal V_\alpha=&\ \mbox{the closure of} \ \{ f\in C^\infty(\mathbb T^n)| \nabla\cdot f=0\} \ \mbox{in} \ \mathcal H^\alpha(\mathbb T^n),\\
\mathcal V'_\alpha=&\ \mbox{the dual of}\ \mathcal V_\alpha.
\end{split}
\end{equation}
The inner product in $L^2(\mathbb T^n)$ is denoted by 
\[\langle f, g \rangle =\int_{\mathbb T^n} f\cdot g \, dx.\]

\medskip

\subsection{Notion of randomization}

We first recollect the large deviation estimates established in \cite{BT2}. 

\begin{Lemma}\label{BT}
Let $(l_i(\omega))_{i=1}^\infty$ be a sequence of real-valued, zero-mean and independent random variables on a probability space $(\Omega, \mathcal A, P)$ with associated distributions $(\mu_i)_{i=1}^\infty$. Assume that there exists $c>0$ such that
\begin{equation}\label{ass1}
\left|\int_{-\infty}^{\infty}e^{\gamma x}\, d\mu_i(x)\right|\leq e^{c\gamma^2} \ \ \ \ \forall \gamma\in\mathbb R \ \ \ \forall i\geq 1.
\end{equation}
Then there exists $\beta>0$ such that 
\begin{equation}\notag
P\left(\omega: \left|\sum_{i=1}^\infty c_il_i(\omega)\right|>\lambda\right)\leq 2e^{-\frac{\beta \lambda^2}{\sum_{i=1}^\infty c_i^2}}   \ \ \ \ \forall \lambda>0 \ \ \ \forall (c_i)_{i=1}^\infty\in \ell^2.
\end{equation}
Consequently, there exists another constant $c>0$ such that 
\begin{equation}\notag
\left\|\sum_{i=1}^\infty c_il_i(\omega)\right\|_{L^q(\Omega)}\leq c \sqrt q\left(\sum_{i=1}^\infty c_i^2\right)^{\frac12}
 \ \ \ \ \forall q\geq 2 \ \ \ \forall (c_i)_{i=1}^\infty\in \ell^2.
\end{equation}
\end{Lemma}

We point out that both the standard Gaussian and Bernoulli variables satisfy the assumption (\ref{ass1}), see \cite{BT2}.

We follow \cite{NPS} to introduce the diagonal randomization on the Sobolev space $\mathcal H^s(\mathbb T^n)$ as follows. 

\begin{Definition}\label{def-rand}
Let $(l_k(\omega))_{k\in \mathbb Z^n}$ be a sequence of real-valued and independent random variables on the probability space $(\Omega,\mathcal A, P)$ as in Lemma \ref{BT}. Let $e_k(x)=e^{ik\cdot x}$ for any $k\in \mathbb Z^n$. For a vector field $f=(f_1, f_2, ..., f_n)\in \mathcal H^s(\mathbb T^n)$ with Fourier coefficients $(a_{k})_{k\in \mathbb Z^n}$ and $a_k=(a_k^1, a_k^2, ..., a_k^n)$, the map 
\begin{equation}\label{R}
\begin{split}
\mathcal R:  (\Omega, \mathcal A) & \longrightarrow \mathcal H^s(\mathbb T^n)\\
\omega & \longrightarrow f^\omega, \ \ \ f^\omega(x)=\left(\sum_{k\in\mathbb Z^n}l_k(\omega)a_k^1 e_k(x), ... , \sum_{k\in\mathbb Z^n}l_k(\omega)a_k^n e_k(x)\right)
\end{split}
\end{equation}
equipped with the Borel sigma algebra is introduced. The map $\mathcal R$ is called randomization. 
\end{Definition}

It follows from Lemma \ref{BT} that the map $\mathcal R$ is measurable and $f^\omega$ is an $\mathcal H^s(\mathbb T^n)$-valued random variable. Moreover, we have 
\begin{equation}\notag
f^\omega\in L^2(\Omega; \mathcal H^s(\mathbb T^n), \ \ \ \|f^\omega\|_{\mathcal H^s}\sim \|f\|_{\mathcal H^s}. 
\end{equation}
Indeed, as shown in \cite{BT2}, the randomization $\mathcal R$ does not provide regularization of $\mathcal H^s$ in term of the regularity index $s$. Nevertheless, it gives rise to improved $L^p$ estimate almost surely.

\medskip

\subsection{Statement of the main results}

\begin{Definition}\label{def-sol1}
Let $\alpha>1$. Let $f\in \mathcal H^{-s}(\mathbb T^n)$ with $s>0$ and 
\begin{equation}\notag
\nabla\cdot f=0, \ \ \ \int_{\mathbb T^n} f\, dx=0.
\end{equation}
A function $B(x,t)$ is said to be a weak solution of the electron MHD (\ref{emhd}) with initial data $f$ on $[0,T]$ if 
\begin{equation}\notag
\langle\frac{d B}{dt}, \phi\rangle+\left<\nabla^\alpha B, \nabla^\alpha \phi\right>+\left<(\nabla\times B)\times B, \nabla\times \phi\right>=0 \ \ \ \mbox{for a.e. } \ \ t \ \ \mbox{and for all } \ \ \phi\in \mathcal V'_\alpha,
\end{equation}
\[B\in L^2_{loc}((0,T); \mathcal V_\alpha (\mathbb T^n))\cap L_{loc}^\infty((0,T); \mathcal H(\mathbb T^n))\cap C_{weak}((0,T); \mathcal H^{-s}(\mathbb T^n)),\]
\[\frac{d B}{dt}\in L^1_{loc}((0,T); \mathcal V'_\alpha (\mathbb T^n)),\]
and
\begin{equation}\notag
\lim_{t\to 0^+}B(t)= f \ \ \mbox{weakly in } \ \ \mathcal H^{-s}(\mathbb T^n).
\end{equation}
\end{Definition}

\begin{Theorem}\label{thm-2d}
Let $f$ be as in Definition \ref{def-sol1}. Let $\alpha\in [\frac43, \frac32]$. Assume $s\in(0, 2\alpha-\frac52)$. There exists a set $\Sigma\subset \Omega$ with $P(\Sigma)=1$ such that for any $\omega\in\Sigma$ the electron MHD (\ref{emhd}) with initial data $f^\omega$ on $\mathbb T^2$ has a global in time weak solution $B$ of the form 
\begin{equation}\label{sol-split}
B=B_{f^\omega}+H
\end{equation}
with $B_{f^\omega}=e^{-t(-\Delta)^\alpha} f^\omega$ and 
\begin{equation}\notag
H\in L^\infty([0,\infty); L^2(\mathbb T^2))\cap L^2([0,\infty); \mathcal H^\alpha(\mathbb T^2)). 
\end{equation}
In addition, if $\alpha\geq \frac32$, the solution is regular and unique. 
\end{Theorem}

\begin{Theorem}\label{thm-3d}
Consider (\ref{emhd}) on $\mathbb T^3$. Let $\alpha\in (\frac{11}{8}, \frac74]$. Assume $s\in (0, 2\alpha-\frac{11}{4})$. Then the first statement of Theorem \ref{thm-2d} holds.
In addition, if $\alpha\geq \frac74$, the solution is regular and unique. 
\end{Theorem}


\bigskip

\section{Outline of the proof of the main results}\label{sec-outline}

The strategy of showing existence of weak solutions for the electron MHD (\ref{emhd}) with initial data $f^\omega$ is to look for solutions in the form $B=B_{f^\omega}+H$ with the linear part 
\[B_{f^\omega}=e^{-t(-\Delta)^\alpha} f^\omega, \ \ \ B_{f^\omega}(x,0)=f^\omega(x)\]
and the remaining nonlinear part $H$.  Denote the bilinear operator
\begin{equation}\notag
\mathcal B(u,v)=(\nabla\times u)\times v.
\end{equation}
Note that 
\[\mathcal B(u,u)= (u\cdot\nabla )u-\nabla \frac{|u|^2}{2}. \]
If $\nabla\cdot u=0$, we can further write 
\[\mathcal B(u,u)= \nabla\cdot(u\otimes u)-\nabla \frac{|u|^2}{2}. \]
One can check that if $B$ satisfies (\ref{emhd}) with initial data $f^\omega$, the nonlinear part $H$ solves the Cauchy problem
\begin{equation}\label{H}
\begin{split}
H_t+\nabla\times\mathcal B(H, H)+\nabla\times\mathcal B(H, B_{f^\omega})\ \ \ &\\
+\nabla\times\mathcal B(B_{f^\omega}, H)+\nabla\times\mathcal B(B_{f^\omega}, B_{f^\omega})=&-(-\Delta)^\alpha H,\\
\nabla\cdot H=&\ 0,\\
H(x, 0)=&\ 0.
\end{split}
\end{equation}
In order to prove Theorems \ref{thm-2d} and \ref{thm-3d}, it is sufficient to show existence of global in time weak solutions for (\ref{H}) on $\mathbb T^n$ with $n=2,3$. We thus proceed to define weak solutions for (\ref{H}) and formulate the existence theorem. 

\begin{Definition}\label{def2}
A function $H(x,t)$ is said to be a weak solution of (\ref{H}) on $[0,T]$ if 
\begin{equation}\notag
\begin{split}
\langle\frac{d H}{dt}, \phi\rangle&+\left<\nabla^\alpha H, \nabla^\alpha \phi\right>+\left<\mathcal B(H, H), \nabla\times \phi\right>\\
&+\left<\mathcal B(H, B_{f^\omega}), \nabla\times \phi\right> +\left<\mathcal B(B_{f^\omega}, H), \nabla\times \phi\right>
+\left<\mathcal B(B_{f^\omega}, B_{f^\omega}), \nabla\times \phi\right>=0 \\
&\ \ \ \mbox{for a.e. } \ \ t \ \ \mbox{and for all } \ \ \phi\in \mathcal V'_\alpha,
\end{split}
\end{equation}
\[H\in L^2((0,T); \mathcal V_\alpha (\mathbb T^n))\cap L^\infty((0,T); \mathcal H(\mathbb T^n)),\ \ \ \ \frac{d H}{dt}\in L^1((0,T); \mathcal V'_\alpha(\mathbb T^n)),\]
and
\begin{equation}\notag
\lim_{t\to 0^+}H(t)= 0 \ \ \mbox{weakly in } \ \ \mathcal H^{-s}(\mathbb T^n).
\end{equation}
\end{Definition}

\medskip

Denote
\begin{equation}\label{B-norms}
\begin{split}
 B_{f^\omega}(\alpha, \beta, s, \gamma, T)
:=&\ \|t^\gamma B_{f^\omega}\|_{L^{p_1}([0,T]; L^{q_1}(\mathbb T^n))}+\|t^\gamma B_{f^\omega}\|_{L^{p_2}([0,T]; L^{q_2}(\mathbb T^n))}\\
&+\|t^\gamma (-\Delta)^{\frac{2-\alpha+\beta}{2}} B_{f^\omega}\|_{L^{p_3}([0,T]; L^{q_3}(\mathbb T^n))}\\
&+\|t^\gamma (-\Delta)^{2-\alpha} B_{f^\omega}\|_{L^{p_4}([0,T]; L^{q_4}(\mathbb T^n))}
\end{split}
\end{equation}
where the parameters $p_i$ and $q_i$ with $1\leq i\leq 4$ are given by
\[p_1=\frac{2\alpha}{2\alpha-2-\epsilon} \ \ \mbox{for a small enough}\ \epsilon>0, \ \ \ \ q_1\gg 1,\]
\[p_2=q_2=\frac{4\alpha}{2\alpha-\beta-2}, \ \ p_3=q_3=\frac{4\alpha}{\beta+2},\ \ p_4=2, \ \ q_4\gg 1.\]

\begin{Theorem}\label{thm-H}
Fix $\lambda>0$. For $n=2$, let $\alpha\in [\frac43, \frac32]$, $s\in(0, 2\alpha-\frac52)$ and $\gamma<0$ such that $0<s<2\alpha-\frac52+2\alpha\gamma$. For $n=3$, let $\alpha\in (\frac{11}{8}, \frac74]$, $s\in(0, 2\alpha-\frac{11}{4})$ and $\gamma<0$ such that $0<s<2\alpha-\frac{11}{4}+2\alpha\gamma$. 
Assume the free evolution $B_{f^\omega}$ satisfies 
\begin{equation}\label{ass-bf1}
\begin{split}
\|B_{f^\omega}\|_{L^2(\mathbb T^n)}\lesssim&\ (1+t^{-\frac{s}{2\alpha}}), \\
\|\nabla^m B_{f^\omega}\|_{L^\infty(\mathbb T^n)}\lesssim&\ (\max\{t^{-\frac12}, t^{-(\frac{2m+n+2s}{2\alpha})}\})^{\frac12} \ \ \mbox{for} \ \ m=0,1, 2
\end{split}
\end{equation}
and 
\begin{equation}\label{ass-bf2}
 B_{f^\omega}(\alpha, \beta, s, \gamma, T) \leq \lambda.
\end{equation}
Then there exists a weak solution $H(x,t)$ to the Cauchy problem (\ref{H}) in the sense of Definition \ref{def2}. The solution is unique in 2D for $\alpha\geq \frac32$ and in 3D for $\alpha\geq \frac74$.
\end{Theorem}

\textbf{Proof of Theorems \ref{thm-2d} and \ref{thm-3d}:} Under the conditions on $\alpha$ and $s$ of Theorems \ref{thm-2d} and \ref{thm-3d}, one can find an appropriate constant $\gamma<0$ such that the assumptions on the parameters of Theorem \ref{thm-H} are satisfied. By Lemmas \ref{le-heat3} and \ref{le-heat4}, the assumption (\ref{ass-bf2}) is satisfied almost surely.  On the other hand, assumption (\ref{ass-bf1}) is guaranteed by Lemma \ref{le-lin1}. Thus the existence of a global weak solution $H(x,t)$ to system (\ref{H}) follows from Theorem \ref{thm-H}. Consequently we obtain the existence of a global weak solution $B(x,t)=B_{f^\omega}(x,t)+H(x,t)$ to (\ref{emhd}) almost surely. Recall that system (\ref{emhd}) is critical in 2D for $\alpha=\frac32$ and in 3D for $\alpha=\frac74$. Hence, above the critical value of $\alpha$, regularity of the solution can be established by standard bootstrapping argument. The proof of uniqueness is presented in Appendix. 

\cbdu

\medskip

The remaining part of the paper is devoted to the proof of Theorem \ref{thm-H}. 
The first step is to establish estimates on the linear part $B_{f^\omega}$ such that assumptions of the theorem are satisfied. This will be the content of Section \ref{sec-linear}. The crucial idea of adapting randomized initial data is revealed in this part. In fact, although the initial data $f$ is merely in $\mathcal H^{-s}$ for $s>0$, the free evolution of the randomized data $f^\omega$ has almost surely improved $L^p$ estimates. As a consequence, we are able to establish suitable a priori estimates for $H$ in Section \ref{sec-est}. 
Then in Section \ref{sec-galerkin} we construct Galerkin approximating solutions for (\ref{H}) by standard arguments, for instance see \cite{CF, DG}, and pass to a limit by applying the a priori estimates.

\bigskip

\section{The linear equation with randomized initial data}\label{sec-linear}

We consider the linear equation with randomized initial data
\begin{equation}\label{eq-lin}
\begin{split}
B_t+(-\Delta)^\alpha B=&\ 0,\\
B(x,0)=&\ f^\omega,
\end{split}
\end{equation}
and establish some a priori estimates for its solution $B_{f^\omega}=e^{-t(-\Delta)^\alpha}f^\omega$. 

The following lemma concerns deterministic estimates.
\begin{Lemma}\label{le-lin1} 
Let $s\geq 0$ and $f^\omega\in \mathcal H^{-s}(\mathbb T^n)$. Then the estimate
\begin{equation}\label{est-lin1}
\|\nabla^m B_{f^\omega}\|_{L^2(\mathbb T^n)}\lesssim \left(1+t^{-\frac{m+s}{2\alpha}}\right)\|f\|_{\mathcal  H^{-s}(\mathbb T^n)},
\end{equation}
holds for any nonnegative integer $m$ and $\alpha>0$; and 
\begin{equation}\label{est-lin2}
\|\nabla^m B_{f^\omega}\|_{L^\infty(\mathbb T^n)}\lesssim \left(\max\{t^{-\frac12}, t^{-\frac{2m+n+2s}{2\alpha}}\}\right)^{\frac12}\|f\|_{\mathcal H^{-s}(\mathbb T^n)}
\end{equation}
holds for $m\geq 0$, $\alpha>0$ and $2m+n\geq \alpha$.
\end{Lemma}
\pf
Note that $y^ae^{-y}\leq C$ for $a\geq 0$ and $y\geq 0$. By Plancherel's theorem we deduce
\begin{equation}\notag
\begin{split}
\|\nabla^m B_{f^\omega}(\cdot, t)\|_{L^2(\mathbb T^n)}\sim&\ \||\xi|^m e^{-|\xi|^{2\alpha}t}\widehat{f^\omega}(\xi)\|_{\ell^2_{\xi}}\\
\sim&\ \|t^{-\frac{m+s}{2\alpha}}(|\xi|^{2\alpha}t)^{\frac{m+s}{2\alpha}} e^{-|\xi|^{2\alpha}t}|\xi|^{-s}\widehat{f^\omega}(\xi)\|_{\ell^2_{\xi}}\\
\lesssim&\ \left(1+t^{-\frac{m+s}{2\alpha}}\right)\|f\|_{\mathcal  H^{-s}(\mathbb T^n)}
\end{split}
\end{equation}
which verifies (\ref{est-lin1}). 

In order to show (\ref{est-lin2}), denote 
\begin{equation}\notag
I= \int_0^\infty (1+\rho^2)^s\rho^{2m}e^{-2\rho^{2\alpha}t} \rho^{n-1}\, d\rho.
\end{equation}
Applying Fourier transform on $\mathbb T^n$ and H\"older's inequality we have 
\begin{equation}\label{est-i0}
\begin{split}
|\nabla^m B_{f^\omega}(x,t)|\leq& \sum_{\xi} |\xi|^m e^{-|\xi|^{2\alpha}t}|\widehat{f^\omega}(\xi)|\\
\lesssim &\ \|f\|_{\mathcal H^{-s}(\mathbb T^n)}\left(\sum_{\xi}(1+|\xi|^2)^{s}|\xi|^{2m} e^{-2|\xi|^{2\alpha}t}\right)^{\frac12}\\
\lesssim &\ \|f\|_{\mathcal H^{-s}(\mathbb T^n)} I^{\frac12}.
\end{split}
\end{equation}
The task now is to estimate the integral $I$.  Changing variable $y=\rho^\alpha\sqrt t$ in the integral we can write
\begin{equation}\notag
\begin{split}
I=&\ \int_0^\infty\frac{1}{\alpha}\left(1+\left(\frac{y^2}{t}\right)^{\frac1\alpha}\right)^s\left(\frac{y}{\sqrt t}\right)^{\frac{2m+n-\alpha}{\alpha}}\frac{1}{\sqrt t}e^{-2y^2}\, dy\\
=& \int_{0}^{\sqrt t} \cdot\cdot\cdot \, dy+\int_{\sqrt t}^\infty \cdot\cdot\cdot \, dy\\
=:&\ I_1+I_2
\end{split}
\end{equation}
where the integrand of $I_1$ and $I_2$ is the same as that of $I$.  For $0\leq y\leq \sqrt t$, we have 
\begin{equation}\notag
1+\left(\frac{y^2}{t}\right)^{\frac1\alpha}\leq C, \ \ \
\left(\frac{y}{\sqrt t}\right)^{\frac{2m+n-\alpha}{\alpha}}\leq 1
\end{equation}
since $2m+n\geq \alpha>0$. Thus the integral $I_1$ satisfies
\begin{equation}\label{est-i1}
\begin{split}
I_1\lesssim & \int_0^{\sqrt t}\frac{1}{\sqrt t}e^{-2y^2}\, dy\\
\lesssim &\ t^{-\frac12} \int_0^{\sqrt t}e^{-2y^2}\, dy\\
\lesssim &\ t^{-\frac12}.
\end{split}
\end{equation}
While for $y>\sqrt t$, it follows 
\begin{equation}\notag
\left(1+\left(\frac{y^2}{t}\right)^{\frac1\alpha}\right)^s\lesssim y^{\frac{2s}{\alpha}}t^{-\frac{s}{\alpha}}
\end{equation}
and hence
\begin{equation}\label{est-i2}
\begin{split}
I_2\lesssim & \int_{\sqrt t}^\infty y^{\frac{2s}{\alpha}}t^{-\frac{s}{\alpha}}\left(\frac{y}{\sqrt t}\right)^{\frac{2m+n-\alpha}{\alpha}}\frac{1}{\sqrt t}e^{-2y^2}\, dy\\
\lesssim &\ t^{-\frac{2m+n+2s}{2\alpha}} \int_{\sqrt t}^\infty y^{\frac{2m+n-\alpha+2s}{\alpha}}e^{-2y^2}\, dy\\
\lesssim &\ t^{-\frac{2m+n+2s}{2\alpha}}.
\end{split}
\end{equation}
Therefore estimate (\ref{est-lin2}) follows from (\ref{est-i0}), (\ref{est-i1}) and (\ref{est-i2}).

\cbdu

\medskip

Probabilistic estimates are obtained as well. Namely,
\begin{Lemma}\label{le-lin2} 
Fix $r\geq p\geq q\geq 2$, $\sigma\geq 0$ and $\gamma\in\mathbb R$ such that $q(\frac{\sigma+s}{2\alpha}-\gamma)<1$. Then for any $T>0$ and $s\geq0$ there exists $C_T(p,q,r,\sigma,\gamma,s)>0$ such that for any $f^\omega\in \mathcal H^{-s}(\mathbb T^n)$
\begin{equation}\label{est-lin3}
\|t^\gamma (-\Delta)^{\frac{\sigma}{2}} B_{f^\omega}\|_{L^r(\Omega; L^q([0,T]; L^p(\mathbb T^n)))}\leq C_T\|f\|_{\mathcal H^{-s}(\mathbb T^n)}.
\end{equation}
Denote 
\begin{equation}\label{est-lin4}
E_{\lambda, T, f, \sigma, p}=\{\omega\in\Omega: \|t^\gamma (-\Delta)^{\frac{\sigma}{2}} B_{f^\omega}\|_{L^q([0,T]; L^p(\mathbb T^n))}\geq \lambda\}.
\end{equation}
Then there exists $c_1>0$ and $c_2>0$ such that
\begin{equation}\label{est-lin5}
P(E_{\lambda, T, f, \sigma, p})\leq c_1 \exp\left\{-\frac{c_2\lambda^2}{C_T\|f\|_{\mathcal H^{-s}}^2}\right\} \ \ \ \forall \lambda>0 \ \ \ \forall f^\omega\in \mathcal H^{-s}(\mathbb T^n).
\end{equation}
\end{Lemma}
\pf
Denote $\langle-\Delta\rangle$ by the operator with Fourier symbol $\widehat{\langle-\Delta\rangle}=1+|\xi|^2$. 
We express the term in Fourier representation 
\begin{equation}\label{est-plin1}
\begin{split}
&t^\gamma (-\Delta)^{\frac{\sigma}{2}} B_{f^\omega}\\
=&\ t^\gamma (-\Delta)^{\frac{\sigma}{2}} \langle-\Delta\rangle^{\frac s2} e^{-t(-\Delta)^\alpha}\langle-\Delta\rangle^{-\frac s2} f^\omega\\
=&\ t^\gamma \sum_{\xi\in \mathbb Z^n}|\xi|^\sigma (1+|\xi|^2)^{\frac{s}2} e^{-t|\xi|^{2\alpha}} (1+|\xi|^2)^{-\frac{s}2}\widehat{f^\omega}(\xi) e_{\xi}(x)\\
\leq &\ t^\gamma \sum_{\xi\in \mathbb Z^n, |\xi|\leq 2}|\xi|^\sigma (1+|\xi|^2)^{\frac{s}2} e^{-t|\xi|^{2\alpha}} (1+|\xi|^2)^{-\frac{s}2}\widehat{f^\omega}(\xi) e_{\xi}(x)\\
&+t^\gamma \sum_{\xi\in \mathbb Z^n, |\xi|> 2}|\xi|^\sigma (1+|\xi|^2)^{\frac{s}2} e^{-t|\xi|^{2\alpha}} (1+|\xi|^2)^{-\frac{s}2}\widehat{f^\omega}(\xi) e_{\xi}(x)\\
=:&\ J_1+J_2.
\end{split}
\end{equation}
Using again the fact that $y^ae^{-y}\leq C$ for $a\geq 0$ and $y>0$, we have
\begin{equation}\notag
\begin{split}
J_1\lesssim&\ t^{\gamma-\frac{\sigma}{2\alpha}} \sum_{\xi\in \mathbb Z^n}(t|\xi|^{2\alpha})^{\frac{\sigma}{2\alpha}} e^{-t|\xi|^{2\alpha}} (1+|\xi|^2)^{-\frac{s}2}\widehat{f^\omega}(\xi) e_{\xi}(x)\\
\lesssim&\ t^{\gamma-\frac{\sigma}{2\alpha}} \sum_{\xi\in \mathbb Z^n}(1+|\xi|^2)^{-\frac{s}2}\widehat{f^\omega}(\xi) e_{\xi}(x)
\end{split}
\end{equation}
and 
\begin{equation}\notag
\begin{split}
J_2\lesssim&\ t^{\gamma-\frac{\sigma+s}{2\alpha}} \sum_{\xi\in \mathbb Z^n}(t|\xi|^{2\alpha})^{\frac{\sigma+s}{2\alpha}} e^{-t|\xi|^{2\alpha}} (1+|\xi|^2)^{-\frac{s}2}\widehat{f^\omega}(\xi) e_{\xi}(x)\\
\lesssim&\ t^{\gamma-\frac{\sigma+s}{2\alpha}} \sum_{\xi\in \mathbb Z^n}(1+|\xi|^2)^{-\frac{s}2}\widehat{f^\omega}(\xi) e_{\xi}(x).
\end{split}
\end{equation}
Denote $h=\langle-\Delta\rangle^{-\frac{s}{2}}f$ and hence $\widehat{h^\omega} (\xi)=(1+|\xi|^2)^{-\frac{s}2}\widehat{f^\omega}(\xi)$ in view of the randomization (\ref{R}). 
We estimate the norm of $J_1$ by applying Minkowski's inequality,
\begin{equation}\notag
\begin{split}
\|J_1\|_{L^r(\Omega; L^q([0,T]; L^p(\mathbb T^n)))}\leq &\ C \| t^{\gamma-\frac{\sigma}{2\alpha}} \sum_{\xi\in \mathbb Z^n}\widehat{h^\omega}(\xi) e_{\xi}(x)\|_{L^r(\Omega; L^q([0,T]; L^p(\mathbb T^n)))}\\
\leq &\ C_r \left\|\left(\sum_{\xi\in \mathbb Z^n}\left|t^{\gamma-\frac{\sigma}{2\alpha}}\widehat{h}(\xi) e_{\xi}(x)\right|^2\right)^{\frac12}\right\|_{L^q([0,T]; L^p(\mathbb T^n))}\\
=&\ C_r \left\|\sum_{\xi\in \mathbb Z^n}\left|t^{\gamma-\frac{\sigma}{2\alpha}}\widehat{h}(\xi) e_{\xi}(x)\right|^2\right\|^{\frac12}_{L^{\frac{q}2}([0,T]; L^{\frac{p}2}(\mathbb T^n))}\\
\leq &\ C_r\left(\int_0^T t^{\frac{q}{2}(2\gamma-\frac{\sigma}{\alpha})}\, dt\right)^{\frac{1}{q}} \left\|\sum_{\xi\in \mathbb Z^n}\left|\widehat{h}(\xi) e_{\xi}(x)\right|^2\right\|^{\frac12}_{L^{\frac{p}2}(\mathbb T^n)}.\\
\end{split}
\end{equation}
We further apply Lemma \ref{BT} to deduce
\begin{equation}\label{est-j1}
\begin{split}
\|J_1\|_{L^r(\Omega; L^q([0,T]; L^p(\mathbb T^n)))}
\leq &\ C_{r,p}\left(\int_0^T t^{\frac{q}{2}(2\gamma-\frac{\sigma}{\alpha})}\, dt\right)^{\frac{1}{q}} \left(\sum_{\xi\in \mathbb Z^n}\left|\widehat{h}(\xi)\right|^4\right)^{\frac14}\\
\leq &\ C_{r,p,q}T^{\gamma-\frac{\sigma}{2\alpha}+\frac1q} \left(\sum_{\xi\in \mathbb Z^n}\left|\widehat{h}(\xi)\right|^2\right)^{\frac12}\\
\leq &\ C_{T, r,p,q, \sigma, \gamma, \alpha}\|f\|_{\mathcal H^{-s}(\mathbb T^n)}
\end{split}
\end{equation}
where we need to require $q(\frac{\sigma}{2\alpha}-\gamma)<1$ for the time integral to be finite. Analogously we have 
\begin{equation}\label{est-j2}
\|J_2\|_{L^r(\Omega; L^q([0,T]; L^p(\mathbb T^n)))}
\leq C_{T, r,p,q}\|f\|_{\mathcal H^{-s}(\mathbb T^n)}
\end{equation}
for $q(\frac{\sigma+s}{2\alpha}-\gamma)<1$. Thus the estimate (\ref{est-lin3}) follows from (\ref{est-plin1}), (\ref{est-j1}) and (\ref{est-j2}).

In the end, the estimate (\ref{est-lin5}) follows from Bienaym\'e-Tchebishev's inequality (see Proposition 4.4 of \cite{BT2}) and Lemma \ref{BT}. 

\cbdu


\medskip

The following maximal regularity result for the free evolution equation is needed to establish energy estimate for $H$ in Section \ref{sec-est}. 
\begin{Lemma}\label{le-heat2} 
Let $T>0$ and $f\in L^2((0,T); L^2(\mathbb T^n))$. Denote \[g(x,t)=\int_0^t e^{-(t-s)(-\Delta)^\alpha}(-\Delta)^\alpha f (x, s)\, ds.\] 
Then we have for any $\alpha>0$
\[\|g\|_{L^2((0,T); L^2(\mathbb T^n))}\lesssim \|f\|_{L^2((0,T); L^2(\mathbb T^n))}.\]
\end{Lemma}
\pf The estimate for $\alpha=1$ is classical, for instance see Theorem 7.3 of \cite{LR}.  We follow the lines of \cite{LR} to prove the estimate for general $\alpha>0$. 

Let $G(x)$ be the kernel function of the operator $e^{-(-\Delta)^\alpha}$,
\[G(x)=(2\pi)^{-\frac{n}{2}}\int_{\mathbb T^n} e^{ix\cdot\xi} e^{-|\xi|^{2\alpha}}\, d\xi\]
and $G(x, t)$ the rescaled function 
\[G(x, t)=t^{-\frac{n}{2\alpha}} G\left(\frac{x}{t^{1/2\alpha}}\right), \ \ \ t>0.\]
We extend $G(x, t)$ to the entire time line by setting $G(x, t)=0$ for $t<0$. We then can write $g(x,t)$ as
\begin{equation}\notag
\begin{split}
g(x,t)=& \int_{-\infty}^{\infty} \int_{\mathbb T^n}\frac{1}{t-s} G(x-y, t-s) f(y,s)\, dy ds \\
=& \left(\frac1{t} G(x,t)\right)* f(x,t)
\end{split}
\end{equation}
where the convolution is in both $x$ and $t$. Thus by Young's inequality we have
\begin{equation}\notag
\|g\|_{L^2((0,T); L^2(\mathbb T^n))}\lesssim \left\|\frac1{t} G(x,t)\right\|_{L^1((0,T); L^1(\mathbb T^n))}\|f\|_{L^2((0,T); L^2(\mathbb T^n))}.
\end{equation}
Note that the Fourier transform of $\frac1{t} G(x,t)$ in both $x$ and $t$ is given by
\begin{equation}\notag
\mathcal F\left(\frac1{t} G\right)(\xi, \tau)=-\int_0^\infty |\xi|^{2\alpha} e^{-t|\xi|^{2\alpha}} e^{-i t\tau}\, dt=-\frac{|\xi|^{2\alpha}}{|\xi|^{2\alpha}+i\tau}.
\end{equation}
We observe that 
\begin{equation}\notag
\left|\mathcal F\left(\frac1{t} G\right)(\xi, \tau)\right|\leq 1
\end{equation}
and hence
\begin{equation}\notag
\left\|\frac1{t} G(x,t)\right\|_{L^1((0,T); L^1(\mathbb T^n))}\leq C.
\end{equation}
It then follows
\[\|g\|_{L^2((0,T); L^2(\mathbb T^n))}\lesssim \|f\|_{L^2((0,T); L^2(\mathbb T^n))}.\]


\cbdu

\medskip

We also need the following estimate. 
\begin{Lemma}\label{le-heat0}
Let $T>0$ and $f\in L^2((0,T); L^2(\mathbb T^n))$. Denote \[g(x,t)=\int_0^t e^{-(t-s)(-\Delta)^\alpha}\nabla^m f (x, s)\, ds.\] 
Then we have for $2\alpha> m$
\[\|g(t)\|_{L^2(\mathbb T^n)}\lesssim \|f\|_{L^2((0,T); L^2(\mathbb T^n))} \ \ \ \ \forall t>0.\]
\end{Lemma}
\pf
For any $\varphi\in L^2(\mathbb T^n)$, using integration by parts and H\"older's inequality we have
\begin{equation}\notag
\begin{split}
\left|\langle g(t), \varphi\rangle \right|=&\ \left| \int_0^t \langle f(s), e^{-(t-s)(-\Delta)^\alpha}\nabla^m \varphi \rangle \, ds\right|\\
\lesssim& \left(\int_0^t\int_{\mathbb T^n} f\, dxds\right)^{\frac12}\left(\int_0^t\int_{\mathbb T^n} \left|e^{-(t-s)(-\Delta)^\alpha}\nabla^m \varphi \right|^2\, dxds\right)^{\frac12}\\
\lesssim& \|f\|_{L^2((0,T); L^2(\mathbb T^n))} \| e^{-t(-\Delta)^\alpha}\nabla^m \varphi \|_{L^2((0,T); L^2(\mathbb T^n))}. 
\end{split}
\end{equation}
In view of Plancherel's theorem, we deduce
\begin{equation}\notag
\begin{split}
\| e^{-t(-\Delta)^\alpha}\nabla^m \varphi \|^2_{ L^2(\mathbb T^n)}=&\ (2\pi)^{-n} \int_{\mathbb T^n}\left| e^{-t|\xi|^{2\alpha}}|\xi|^{m}\widehat\varphi^2(\xi)\right|\, d\xi\\
\lesssim &\ \frac{1}{t^{\frac{m}{2\alpha}}} \|\varphi\|_{L^2(\mathbb T^n)}
\end{split}
\end{equation}
where we used the fact $x^ae^{-x^2}\leq C$ for $x>0$ and $a>0$. Therefore, we obtain for $m<2\alpha$
\begin{equation}\notag
\begin{split}
\| e^{-t(-\Delta)^\alpha}\nabla^m \varphi \|^2_{L^2((0,T); L^2(\mathbb T^n))}
\lesssim &\ \|\varphi\|_{L^2(\mathbb T^n)}\int_0^T\frac{1}{t^{\frac{m}{2\alpha}}} \, ds \lesssim  \|\varphi\|_{L^2(\mathbb T^n)}.
\end{split}
\end{equation}
Therefore we have
\begin{equation}\notag
\left|\langle g(t), \varphi\rangle \right|
\lesssim \|f\|_{L^2((0,T); L^2(\mathbb T^n))} \|\varphi\|_{L^2(\mathbb T^n)} \ \ \forall \varphi\in L^2(\mathbb T^n)
\end{equation}
which concludes the proof of the lemma.

\cbdu

\medskip

We introduce one more probabilistic estimate for the free evolution $B_{f^\omega}$ in each case of 2D and 3D. 

\begin{Lemma}\label{le-heat3}
Let $n=2$, $\alpha\in [\frac43, \frac32]$ and $\beta=3-2\alpha$. Let $0<s<2\alpha-\frac52+2\alpha\gamma$ for some $\gamma<0$ such that $2\alpha-\frac52+2\alpha\gamma>0$. Let $ B_{f^\omega}(\alpha, \beta, s, \gamma, T)$ be the sum of the norms defined in (\ref{B-norms}).
There exists a set $\Sigma\subset \Omega$ with $P(\Sigma)=1$ such that for any $\omega\in\Sigma$ we can find a constant $\lambda>0$ such that 
\begin{equation}\notag
 B_{f^\omega}(\alpha, \beta, s, \gamma, T)\leq \lambda.
\end{equation}
\end{Lemma}
\pf
For any $\lambda>0$ denote
\begin{equation}\notag
E(\lambda):=E(\lambda, s, \alpha, f, \gamma, T)=\left\{\omega\in\Omega | B_{f^\omega}(\alpha, \beta, s, \gamma, T)>\lambda \right\}.
\end{equation}
For any $j\geq 0$ we also denote $\lambda_j=2^j$ and $E_j=E(\lambda_j)$. Note that $E_{j+1}\subset E_j$. Take 
\[\Sigma=\cup_{j\geq 0} E_j^{c}\subset \Omega.\]
One can check that the parameters satisfy the condition of Lemma \ref{le-lin2}. Hence it follows from Lemma \ref{le-lin2} that
\begin{equation}\notag
P(E_j)\leq c_1 \exp\left\{-\frac{c_2\lambda_j^2}{C_T\|f\|_{H^{-s}}^2}\right\} \ \ \ \forall j\geq0 \ \ \ \forall f^\omega\in (H^{-s}(\mathbb T^2))^2.
\end{equation}
Therefore we deduce
\begin{equation}\notag
\begin{split}
1\geq &\ P(\Sigma)=1-P(\Sigma^c)=1-P(\cap E_j)=1-P(\lim_{j\to\infty} E_j)\\
\geq &\ 1-\lim_{j\to\infty} c_1 \exp\left\{-\frac{c_2\lambda_j^2}{C_T\|f\|_{H^{-s}}^2}\right\}=1
\end{split}
\end{equation}
which immediately gives $P(\Sigma)=1$. By definition of $\Sigma$, we see that for any $\omega\in\Sigma$ there exists $j\geq 0$ such that $\omega\in E_j^c$, i.e. 
\[B_{f^\omega}(\alpha, \beta, s, \gamma, T)\leq \lambda_j.\]

\cbdu

\medskip

\begin{Lemma}\label{le-heat4}
Let $n=3$, $\alpha\in (\frac{11}{8}, \frac{7}{4}]$ and $\beta=\frac72-2\alpha$. Let $0<s<2\alpha-\frac{11}{4}+2\alpha\gamma$ for some $\gamma<0$ such that $2\alpha-\frac{11}{4}+2\alpha\gamma>0$. 
There exists a set $\Sigma\subset \Omega$ with $P(\Sigma)=1$ such that for any $\omega\in\Sigma$ we can find a constant $\lambda>0$ such that 
\begin{equation}\notag
B_{f^\omega}(\alpha, \beta, s, \gamma, T)\leq \lambda.
\end{equation}
\end{Lemma}

\pf
We observe that the parameters specified in the lemma satisfy the assumptions of Lemma \ref{le-lin2}. The proof follows from an analogous argument as that of Lemma \ref{le-heat3}.

\cbdu

\bigskip

\section{A priori estimates for $H$}\label{sec-est}

In this section we establish a priori estimates for the nonlinear part $H$ which solves the Cauchy problem (\ref{H}). Notice that $B_{f^\omega}$ appears in the quadratic nonlinear terms of (\ref{H}) and the estimates of $B_{f^\omega}$ in (\ref{est-lin1}) and (\ref{est-lin2}) exhibit a singularity at $t=0$. To avoid this singularity, we choose to perform the estimates near time zero by working with the integral form of (\ref{H}). Away from time zero, the estimates can be obtained from (\ref{H}). Therefore, before starting the estimates we introduce the mild formulation of (\ref{H}) and show that the two formulations are equivalent under appropriate assumptions.

Denote 
\begin{equation}\notag
\begin{split}
\tilde Q(x,t)=&\ \nabla\times\left[\mathcal B(H+B_{f^\omega}, H+B_{f^\omega})\right],\\
Q(x,t)=&\ (H+B_{f^\omega})\otimes (H+B_{f^\omega}) (x, t).
\end{split}
\end{equation}
Since $\nabla \cdot H=0$ and $\nabla\cdot B_{f^\omega}=0$, we have that following several vector identities 
\[\tilde Q(x,t)= \nabla\times\nabla\cdot Q(x,t).\]
Thus we can write
\begin{equation}\label{H2}
\begin{split}
H(x,t)=&-\int_0^t e^{-(t-s)(-\Delta)^\alpha} \tilde Q(x, s)\, ds\\
=&-\int_0^t e^{-(t-s)(-\Delta)^\alpha} \nabla\times\nabla\cdot Q(x, s)\, ds.
\end{split}
\end{equation}

\begin{Lemma}\label{le-equiv}
Assume $B_{f^\omega}$ satisfies the assumptions (\ref{ass-bf1}) and (\ref{ass-bf2}). Then $H$ is a weak solution to (\ref{H}) if and only if $H\in L^\infty((0,T); \mathcal H(\mathbb T^n))\cap L^2((0,T); \mathcal V_\alpha (\mathbb T^n))$ is a solution to (\ref{H2}). 
\end{Lemma}
\pf
We follow the lines of \cite{LR}. First we assume $H\in L^\infty((0,T); \mathcal H(\mathbb T^n))\cap L^2((0,T); \mathcal V_\alpha (\mathbb T^n))$ is a solution to (\ref{H2}).  Denote 
\begin{equation}\label{eq-M}
\mathcal M (H)(x,t)=-\int_0^t e^{-(t-s)(-\Delta)^\alpha} \nabla\times\nabla\cdot Q(x, s)\, ds.
\end{equation}
Thanks to the assumptions (\ref{ass-bf1}) and (\ref{ass-bf2}) and the fact $H\in L^\infty((0,T); \mathcal H(\mathbb T^n))\cap L^2((0,T); \mathcal V_\alpha (\mathbb T^n))$ we have $Q\in L^1((0,T); L^1(\mathbb T^n))$ and hence 
\[\nabla\times\nabla\cdot Q \in L^1((0,T); \mathcal D').\]
It then follows 
\begin{equation}\notag
e^{-(t-s)(-\Delta)^\alpha} \nabla\times\nabla\cdot Q\in L^1((0,T); C^\infty(\mathbb T^n)).
\end{equation}
Thus by Leibniz rule we have
\begin{equation}\notag
\partial_t \mathcal M (H)(x,t)=-(-\Delta)^\alpha \mathcal M (H)(x,t)- \nabla\times\nabla\cdot Q
\end{equation}
in the distributional sense. On the other hand, we see 
\[\lim_{t\to 0^+} H(x,t)=0. \]
Therefore $H=\mathcal M(H)$ is a weak solution of (\ref{H}).

Conversely, we assume $H$ is a weak solution of (\ref{H}). Define $\mathcal M(x,t)$ as in (\ref{eq-M}). Applying the estimates from Proposition \ref{prop} below near time zero we obtain
\[\mathcal M\in L^\infty((0,T); \mathcal H(\mathbb T^n))\cap L^2((0,T); \mathcal V_\alpha(\mathbb T^n)),\]
\[\frac{d\mathcal M}{dt}\in L^1((0,T); \mathcal V'_\alpha(\mathbb T^n)).\]
Hence we deduce by Leibniz rule again
\begin{equation}\notag
\begin{split}
\partial_t(\mathcal M-H)=&-(-\Delta)^\alpha \mathcal M (H)- \nabla\times\nabla\cdot Q+(-\Delta)^\alpha \mathcal H+ \nabla\times\nabla\cdot Q\\
=&-(-\Delta)^\alpha (\mathcal M (H)-H)
\end{split}
\end{equation}
which is satisfied in the distributional sense. Note that 
\[\lim_{t\to 0^+}(\mathcal M (H)(t)-H(t))=0.\]
It then follows from the uniqueness of the generalized heat flow that $\mathcal M=H$ and hence $H$ is a weak solution of (\ref{H2}).

\cbdu

\medskip

Denote the basic energy functional 
\begin{equation}\notag
\mathcal E(H)(t)=\|H(t)\|_{L^2(\mathbb T^n)}^2+2\int_0^t \int_{\mathbb T^n} |\nabla^\alpha H(s)|^2\, dx\, ds,
\end{equation}
and the higher oder energy functional for some $\beta$ to be determined
\begin{equation}\label{E2}
\begin{split}
\mathcal E_{\alpha}(H)(t)=&\ \mathcal E(H)(t)+\mathcal E((-\Delta)^{\frac{\beta}{2}}H)(t)\\
=&\ \|H(t)\|_{L^2(\mathbb T^n)}^2+2\int_0^t \int_{\mathbb T^n} |\nabla^\alpha H(s)|^2\, dx\, ds\\
&+\|H(t)\|_{\mathcal H^{\beta}(\mathbb T^n)}^2+2\int_0^t \int_{\mathbb T^n} |\nabla^{\alpha+\beta} H(s)|^2\, dx\, ds.
\end{split}
\end{equation}

\begin{Proposition}\label{prop}
Assume $B_{f^\omega}$ satisfies the conditions (\ref{ass-bf1}) and (\ref{ass-bf2}). Let $H\in L^\infty((0,T); \mathcal H(\mathbb T^n))\cap L^2((0,T); \mathcal V_\alpha(\mathbb T^n))$ be a solution to (\ref{H}). Then there exists a constant $C(T, \lambda, s)$ such that
\begin{equation}\label{ap-est1} 
\mathcal E(H)(t)\leq C(T, \lambda, s) \ \ \ \mbox{for all} \ \ t\in[0,T], 
\end{equation}
and 
\begin{equation}\label{ap-est2} 
\begin{split}
\left\|\frac{d}{dt}H\right\|_{L^2((0,T); \mathcal H^{-2\alpha}(\mathbb T^2))}\leq &\ C(T, \lambda, s), \\
\left\|\frac{d}{dt}H\right\|_{L^{\frac{4\alpha}3}((0,T); \mathcal H^{-2\alpha}(\mathbb T^3))}\leq &\ C(T, \lambda, s).
\end{split}
\end{equation}
\end{Proposition}
\pf
As discussed earlier, in order to obtain the estimate (\ref{ap-est1}) we split the time interval into two regimes $[0, t_0]$ and $[t_0, T]$ for a small time $t_0>0$ to be determined later. On $[0,t_0]$ we work with the integral form (\ref{H2}) and take the advantage of the fact $H(x, 0)=0$; while on $[t_0, T]$ we work with the differential form (\ref{H}) since no time singularity presents on this interval. Achieving the estimates on $[0,t_0]$ turns out to be more challenging. We apply the higher order energy method to overcome the obstruction by estimating the energy functional $\mathcal E_{\alpha}$ instead of $\mathcal E$. We choose to treat the 2D and 3D cases separately. Thus the proof consists four parts: (i) estimate of $\mathcal E_\alpha$ on $[0,t_0]$ in 2D; (ii) estimate of $\mathcal E_\alpha$ on $[0,t_0]$ in 3D; (iii) estimate of $\mathcal E$ on $[t_0, T]$ for any spatial dimension; (iv) estimate of $\frac{d}{dt}H$. 

\medskip

{\textbf{(i) Estimates on $[0,t_0]$ in 2D.}}
By Lemma \ref{le-heat0} we have for $\alpha>1$ and any $0<t\leq t_0$
\begin{equation}\label{L22}
\begin{split}
\|H(t)\|_{L^2(\mathbb T^2)}\lesssim &\ \|Q\|_{L^2((0, t_0); L^{2}(\mathbb T^2))}, \\
\|H(t)\|_{\mathcal H^{\beta}(\mathbb T^2)}\lesssim &\ \|Q\|_{L^2((0, t_0); \mathcal H^{\beta}(\mathbb T^2))}.  \\
\end{split}
\end{equation}
In view of the second line of (\ref{H2}) we have 
\begin{equation}\notag
\begin{split}
(-\Delta)^{\frac{\alpha}{2}}H(x,t)=&-\int_0^t e^{-(t-s)(-\Delta)^\alpha}(-\Delta)^{\alpha} \nabla\times\nabla\cdot (-\Delta)^{-\frac{\alpha}{2}}Q(x, s)\, ds\\
(-\Delta)^{\frac{\alpha+\beta}{2}}H(x,t)=&-\int_0^t e^{-(t-s)(-\Delta)^\alpha}(-\Delta)^{\alpha} \nabla\times\nabla\cdot (-\Delta)^{\frac{\beta-\alpha}{2}}Q(x, s)\, ds\\
\end{split}
\end{equation}
and hence we have  from Lemma \ref{le-heat2}
\begin{equation}\label{Lpq2}
\begin{split}
\|H(t)\|_{L^2((0,t_0); \mathcal H^\alpha(\mathbb T^2))} \lesssim&\ \|Q\|_{L^2((0, t_0); \mathcal H^{2-\alpha}(\mathbb T^2))},\\
\|H(t)\|_{L^2((0, t_0); \mathcal H^{\alpha+\beta}(\mathbb T^2))} \lesssim&\ \|Q\|_{L^2((0, t_0); \mathcal H^{2-\alpha+\beta}(\mathbb T^2))}.
\end{split}
\end{equation}
In view of (\ref{L22}) and (\ref{Lpq2}), we need to estimate $\|Q\|_{L^2((0, t_0); L^{2}(\mathbb T^2))}$, $\|Q\|_{L^2((0, t_0); \mathcal H^{2-\alpha}(\mathbb T^2))}$, $\|Q\|_{L^2((0, t_0); \mathcal H^{\beta}(\mathbb T^2))}$, and $\|Q\|_{L^2((0, t_0); \mathcal H^{2-\alpha+\beta}(\mathbb T^2))}$. With the restriction of $1<\alpha\leq \frac32$ and $\beta\geq0$, we have $2-\alpha\leq 2-\alpha+\beta$ and $\beta\leq 2-\alpha+\beta$. Thus it is sufficient to estimate the last one, i.e. $\|Q\|_{L^2((0, t_0); \mathcal H^{2-\alpha+\beta}(\mathbb T^2))}$.

Note that
\begin{equation}\label{est-q1}
\begin{split}
&\|Q\|_{L^2((0, t_0); \mathcal H^{2-\alpha+\beta}(\mathbb T^2))}\\
\lesssim&\ \|H\otimes H\|_{L^2((0, t_0); \mathcal H^{2-\alpha+\beta}(\mathbb T^2))}+\|H\otimes B_{f^\omega}\|_{L^2((0, t_0); \mathcal H^{2-\alpha+\beta}(\mathbb T^2))}\\
&+\|B_{f^\omega}\otimes B_{f^\omega}\|_{L^2((0, t_0); \mathcal H^{2-\alpha+\beta}(\mathbb T^2))}.
\end{split}
\end{equation}
It follows from H\"older's inequality that if $\alpha-\beta\geq 1$ 
\begin{equation}\notag
\begin{split}
&\|H\otimes H\|_{L^2((0, t_0); \mathcal H^{2-\alpha+\beta}(\mathbb T^2))}\\
\lesssim &\ \|H\nabla^{2-\alpha+\beta}H\|_{L^2([0,t_0];L^{2}(\mathbb T^2))}  \\
\lesssim &\ \|H\|_{L^4([0,t_0];L^{\frac{4}{\alpha-\beta-1}}(\mathbb T^2))} \|\nabla^{2-\alpha+\beta} H\|_{L^4([0,t_0];L^{\frac{4}{3-\alpha+\beta}}(\mathbb T^2))}.
\end{split}
\end{equation}
Since by Sobolev embedding 
\[\|H\|_{L^4([0,t_0];L^{\frac{4}{\alpha-\beta-1}}(\mathbb T^2))}\lesssim \|\nabla^{2-\alpha+\beta} H\|_{L^4([0,t_0];L^{\frac{4}{3-\alpha+\beta}}(\mathbb T^2))},\]
we only estimate the latter for $1\leq \alpha-\beta\leq 3$:
\begin{equation}\notag
\begin{split}
&\|\nabla^{2-\alpha+\beta} H\|_{L^4([0,t_0];L^{\frac{4}{3-\alpha+\beta}}(\mathbb T^2))}\\
= & \left(\int_0^{t_0}\|\nabla^{2-\alpha+\beta}H\|^4_{L^{\frac{4}{3-\alpha+\beta}}(\mathbb T^2)} \, dt\right)^{\frac14}\\
\lesssim & \left( \int_0^{t_0}\|\nabla^{\beta}H\|^2_{L^{2}(\mathbb T^2)}\|\nabla^{m+\beta}H\|^2_{L^{2}(\mathbb T^2)} \, dt \right)^{\frac14}\\
\lesssim &\ \left(\sup_{t\in(0, t_0)}\|\nabla^{\beta} H(t)\|^2_{L^2_x}\right)^{\frac14}\left(\int_0^{t_0} \|\nabla^{\alpha+\beta}H\|_{L^2_x}^2\, dt\right)^{\frac14}\\
\lesssim &\ \mathcal E_{\alpha}^{\frac12}(H)(t_0)
\end{split}
\end{equation}
with $m=3-\alpha-\beta\leq \alpha$ provided $\beta\geq 3-2\alpha$. 
Hence, if $1\leq \alpha-\beta\leq 3$ and $\beta\geq 3-2\alpha$ we have
\begin{equation}\label{est-q2}
\|H\otimes H\|_{L^2((0, t_0); H^{2-\alpha+\beta}(\mathbb T^2))}
\lesssim \mathcal E_{\alpha}(H)(t_0).
\end{equation}
The conditions $1\leq \alpha-\beta\leq 3$ and $\beta\geq 3-2\alpha$ imply 
\[\frac43\leq \alpha\leq \frac32.\] 
To optimize the final result, we take the smallest $\beta=3-2\alpha$ from now on.

We continue to estimate
\begin{equation}\notag
\begin{split}
&\|H\otimes B_{f^\omega}\|_{L^2((0, t_0); H^{2-\alpha+\beta}(\mathbb T^2))}\\
\lesssim&\ \|H \nabla^{2-\alpha+\beta} B_{f^\omega}\|_{L^2((0, t_0); L^{2}(\mathbb T^2))}+\|B_{f^\omega} \nabla^{2-\alpha+\beta}H\|_{L^2((0, t_0); L^{2}(\mathbb T^2))}.
\end{split}
\end{equation}
The first term is estimated as follows by applying H\"older's inequality and (\ref{ass-bf2})
\begin{equation}\notag
\begin{split}
&\|H \nabla^{2-\alpha+\beta} B_{f^\omega}\|_{L^2((0, t_0); L^{2}(\mathbb T^2))}\\
\lesssim&\ \|H\|_{L^\infty((0, t_0); L^{p}(\mathbb T^2))}\|\nabla^{2-\alpha+\beta} B_{f^\omega}\|_{L^2((0, t_0); L^{p'}(\mathbb T^2))}\\
\lesssim&\ \|H\|_{L^\infty((0, t_0); H^{\beta}(\mathbb T^2))}\|\nabla^{2-\alpha+\beta} B_{f^\omega}\|_{L^2((0, t_0); L^{p'}(\mathbb T^2))}\\
\lesssim&\ \lambda t_0^{-\gamma}\mathcal E_{\alpha}^{\frac12}(H)(t_0)
\end{split}
\end{equation}
with $\frac{1}{p}+\frac{1}{p'}=\frac12$ and $p=2+\epsilon$ such that the Sobolev embedding holds. 
The second term is estimated as
\begin{equation}\notag
\begin{split}
&\|B_{f^\omega} \nabla^{2-\alpha+\beta} H\|_{L^2((0, t_0); L^{2}(\mathbb T^2))}\\
\lesssim&\ \|\nabla^{2-\alpha+\beta} H\|_{L^{p}((0, t_0); L^{q}(\mathbb T^2))}\|B_{f^\omega}\|_{L^{p'}((0, t_0); L^{q'}(\mathbb T^2))}
\end{split}
\end{equation}
with $\frac{1}{p}+\frac{1}{p'}=\frac12$, $\frac{1}{p}+\frac{1}{p'}=\frac12$, $p'\leq q'$ and $p\geq q$.
By Gagliardo-Nirenberg's inequality we know 
\begin{equation}\notag
\|\nabla^{2-\alpha+\beta} H\|_{L^q(\mathbb T^2)}\lesssim \|\nabla^{\alpha+\beta} H\|_{L^2(\mathbb T^2)}^{\theta}\|\nabla^{\beta} H\|_{L^2(\mathbb T^2)}^{1-\theta}
\end{equation}
with $q=\frac{2}{3-\alpha-\alpha\theta}$ and $p\theta=2$. Take $q=2+\epsilon$ for some small constant $\epsilon>0$, we obtain $p'=\frac{2\alpha}{2\alpha-2-\epsilon}$ for another small constant $\epsilon>0$, and analogous computation as before shows
\begin{equation}\notag
\begin{split}
&\|\nabla^{2-\alpha+\beta} H\|_{L^{p}((0, t_0); L^{q}(\mathbb T^2))}\\
\lesssim &\ \left(\int_0^{t_0} \|\nabla^{\alpha+\beta} H\|_{L^2(\mathbb T^2)}^2 \|\nabla^{\beta}H\|_{L^2(\mathbb T^2)}^{p-2} \, dt\right)^{\frac1p}\\
\lesssim &\ \left(\sup_{0\leq t\leq t_0}\|\nabla^{\beta}H\|_{L^2(\mathbb T^2)}^2 \right)^{\frac{p-2}{2p}}
\left(\int_0^{t_0} \|\nabla^{\alpha+\beta} H\|_{L^2(\mathbb T^2)}^2 \, dt\right)^{\frac{1}{p}}\\
\lesssim&\ \mathcal E_{\alpha}^{\frac12}(H)(t_0).
\end{split}
\end{equation}
Hence we have
\begin{equation}\notag
\|B_{f^\omega} \nabla^{2-\alpha+\beta} H\|_{L^2((0, t_0); L^{2}(\mathbb T^2))}
\lesssim \lambda t_0^{-\gamma}\mathcal E_{\alpha}^{\frac12}(H)(t_0).
\end{equation}
Collecting the estimates above we obtain
\begin{equation}\label{est-q3}
\|H\otimes B_{f^\omega}\|_{L^2((0, t_0); H^{2-\alpha+\beta}(\mathbb T^2))}
\lesssim \lambda t_0^{-\gamma}\mathcal E_{\alpha}^{\frac12}(H)(t_0).
\end{equation}

In the end the condition (\ref{ass-bf2}) again implies
\begin{equation}\label{est-q4}
\begin{split}
&\|B_{f^\omega}\otimes B_{f^\omega}\|_{L^2((0, t_0); H^{2-\alpha+\beta}(\mathbb T^2))}\\
\lesssim&\ \|B_{f^\omega} \nabla^{2-\alpha+\beta}B_{f^\omega}\|_{L^2((0, t_0); L^2(\mathbb T^2))}\\
\lesssim&\ \|B_{f^\omega}\|_{L^p((0, t_0); L^{p}(\mathbb T^2))}\|\nabla^{2-\alpha+\beta} B_{f^\omega}\|_{L^{p'}((0, t_0); L^{p'}(\mathbb T^2))}\\
\lesssim&\ \lambda^2t_0^{-2\gamma}
\end{split}
\end{equation}
with $\frac{1}{p}+\frac{1}{p'}=\frac12$.  This estimate needs to be optimized such that 
\[ p'\left(\frac{2-\alpha+\beta+s}{2\alpha}-\gamma\right)<1, \ \ p\left(\frac{s}{2\alpha}-\gamma\right)<1\] 
for the largest possible value of $s$ and some $\gamma<0$. The optimization results in 
\[p=\frac{4\alpha}{2\alpha-\beta-2}, \ \ \ p'=\frac{4\alpha}{\beta+2}, \ \ \ s<\frac{2\alpha}{p}+2\alpha\gamma=\frac12(2\alpha-\beta-2)+2\alpha\gamma.\]

Combining (\ref{E2}), (\ref{L22}), (\ref{Lpq2}) and the estimates (\ref{est-q1})-(\ref{est-q4}) we obtain for $t\in[0,t_0]$, some $\gamma<0$ and $\alpha, \beta$ and $s$ satisfying
\[\alpha\geq \frac43, \ \ \ \beta=3-2\alpha, \ \ \ 0<s<2\alpha-\frac52+2\alpha\gamma\]
that
\begin{equation}\notag
\begin{split}
\mathcal E^{\frac12}_{\alpha}(H)(t_0)\lesssim&\ \|H\|_{L^\infty((0,t_0); L^2(\mathbb T^2))}+\|H\|_{L^\infty((0,t_0); \mathcal H^{2-\alpha}(\mathbb T^2))}\\
&+\|H\|_{L^2((0,t_0); \mathcal H^\alpha(\mathbb T^2))}+\|H\|_{L^2((0,t_0); \mathcal H^{2-\alpha+\beta}(\mathbb T^2))}\\
\lesssim&\ \mathcal E_{\alpha}(H)(t_0)+\lambda t_0^{-\gamma} \mathcal E^{\frac12}_{\alpha}(H)(t_0)+\lambda^2 t_0^{-2\gamma}\\
\leq&\ C_1\mathcal E_{\alpha}(H)(t_0)+C_2\lambda t_0^{-\gamma} \mathcal E^{\frac12}_{\alpha}(H)(t_0)+C_3\lambda^2 t_0^{-2\gamma} 
\end{split}
\end{equation}
for some constants $C_1, C_2$ and $C_3$. 
By a continuity argument we conclude that for small enough $t_0$ such that $C_3\lambda^2 t_0^{-2\gamma} \ll1$,
\begin{equation}\notag
\mathcal E(H)(t)\leq \mathcal E_{\alpha}(H)(t)\leq C \ \ \ \forall t\in[0,t_0].
\end{equation}

{\textbf{(ii) Estimates on $[0,t_0]$ in 3D.}}
The estimates will be carried on analogously as in the 2D case. Differences come in when we apply Gagliardo-Nirenberg's interpolation inequality and Sobolev embedding inequality. It is again sufficient to estimate
\begin{equation}\label{est-q6}
\begin{split}
&\|Q\|_{L^2((0, t_0); \mathcal H^{2-\alpha+\beta}(\mathbb T^3))}\\
\lesssim&\ \|H\otimes H\|_{L^2((0, t_0); \mathcal H^{2-\alpha+\beta}(\mathbb T^3))}+\|H\otimes B_{f^\omega}\|_{L^2((0, t_0); \mathcal H^{2-\alpha+\beta}(\mathbb T^3))}\\
&+\|B_{f^\omega}\otimes B_{f^\omega}\|_{L^2((0, t_0); \mathcal H^{2-\alpha+\beta}(\mathbb T^3))}.
\end{split}
\end{equation}
By H\"older's inequality the first term on the right hand side of (\ref{est-q6}) is estimated for $\alpha-\beta\geq 1/2$ 
\begin{equation}\notag
\begin{split}
&\|H\otimes H\|_{L^2((0, t_0); \mathcal H^{2-\alpha+\beta}(\mathbb T^3))}\\
\lesssim &\ \|H\nabla^{2-\alpha+\beta}H\|_{L^2([0,t_0];L^{2}(\mathbb T^3))}  \\
\lesssim &\ \|H\|_{L^4([0,t_0];L^{\frac{12}{2\alpha-2\beta-1}}(\mathbb T^3))} \|\nabla^{2-\alpha+\beta} H\|_{L^4([0,t_0];L^{\frac{12}{7-2\alpha+2\beta}}(\mathbb T^3))}.
\end{split}
\end{equation}
In view of Sobolev embedding 
\[\|H\|_{L^4([0,t_0];L^{\frac{12}{2\alpha-2\beta-1}}(\mathbb T^3))}\lesssim \|\nabla^{2-\alpha+\beta} H\|_{L^4([0,t_0];L^{\frac{12}{7-2\alpha+2\beta}}(\mathbb T^3))},\]
we only need to estimate for $ \alpha-\beta\geq 1/2$
\begin{equation}\notag
\begin{split}
&\|\nabla^{2-\alpha+\beta} H\|_{L^4([0,t_0];L^{\frac{12}{7-2\alpha+2\beta}}(\mathbb T^3))}\\
= & \left(\int_0^{t_0}\|\nabla^{2-\alpha+\beta}H\|^4_{L^{\frac{12}{7-2\alpha+2\beta}}(\mathbb T^3)} \, dt\right)^{\frac14}\\
\lesssim & \left( \int_0^{t_0}\|\nabla^{\beta}H\|^2_{L^{2}(\mathbb T^3)}\|\nabla^{m+\beta}H\|^2_{L^{2}(\mathbb T^3)} \, dt \right)^{\frac14}\\
\lesssim &\ \left(\sup_{t\in(0, t_0)}\|\nabla^{\beta} H(t)\|^2_{L^2_x}\right)^{\frac14}\left(\int_0^{t_0} \|\nabla^{\alpha+\beta}H\|_{L^2_x}^2\, dt\right)^{\frac14}\\
\lesssim &\ \mathcal E_{\alpha}^{\frac12}(H)(t_0)
\end{split}
\end{equation}
with $m=\frac72-\alpha-\beta\leq \alpha$ provided $\beta\geq \frac72-2\alpha$. 
Hence for $ \alpha-\beta\geq 1/2$ and $\beta\geq \frac72-2\alpha$ (which imply $\alpha\geq 4/3$) we have
\begin{equation}\label{est-q7}
\|H\otimes H\|_{L^2((0, t_0); \mathcal H^{2-\alpha+\beta}(\mathbb T^3))}
\lesssim \mathcal E_{\alpha}(H)(t_0).
\end{equation}
As before, we choose the smallest $\beta= \frac72-2\alpha$ from now on.

Following the inequality 
\begin{equation}\notag
\begin{split}
&\|H\otimes B_{f^\omega}\|_{L^2((0, t_0); \mathcal H^{2-\alpha+\beta}(\mathbb T^3))}\\
\lesssim&\ \|H \nabla^{2-\alpha+\beta} B_{f^\omega}\|_{L^2((0, t_0); L^{2}(\mathbb T^3))}+\|B_{f^\omega} \nabla^{2-\alpha+\beta}H\|_{L^2((0, t_0); L^{2}(\mathbb T^3))},
\end{split}
\end{equation}
we proceed to estimate the former one on the right hand side as
\begin{equation}\notag
\begin{split}
&\|H \nabla^{2-\alpha+\beta} B_{f^\omega}\|_{L^2((0, t_0); L^{2}(\mathbb T^3))}\\
\lesssim&\ \|H\|_{L^\infty((0, t_0); L^{p}(\mathbb T^2))}\|\nabla^{2-\alpha+\beta} B_{f^\omega}\|_{L^2((0, t_0); L^{p'}(\mathbb T^3))}\\
\lesssim&\ \|H\|_{L^\infty((0, t_0); H^{\beta}(\mathbb T^3))}\|\nabla^{2-\alpha+\beta} B_{f^\omega}\|_{L^2((0, t_0); L^{p'}(\mathbb T^3))}\\
\lesssim&\ \lambda t_0^{-\gamma}\mathcal E_{\alpha}^{\frac12}(H)(t_0)
\end{split}
\end{equation}
with $\frac{1}{p}+\frac{1}{p'}=\frac12$ and $p=2+\epsilon$ such that the Sobolev embedding holds. 
The latter one is estimated as
\begin{equation}\notag
\begin{split}
&\|B_{f^\omega} \nabla^{2-\alpha+\beta} H\|_{L^2((0, t_0); L^{2}(\mathbb T^3))}\\
\lesssim&\ \|\nabla^{2-\alpha+\beta} H\|_{L^{p}((0, t_0); L^{q}(\mathbb T^3))}\|B_{f^\omega}\|_{L^{p'}((0, t_0); L^{q'}(\mathbb T^3))}
\end{split}
\end{equation}
with $\frac{1}{p}+\frac{1}{p'}=\frac12$ and $\frac{1}{p}+\frac{1}{p'}=\frac12$, $p'\leq q'$ and $p\geq q$.
We use Gagliardo-Nirenberg's inequality 
\begin{equation}\notag
\|\nabla^{2-\alpha+\beta} H\|_{L^q(\mathbb T^3)}\lesssim \|\nabla^{\alpha+\beta} H\|_{L^2(\mathbb T^3)}^{\theta}\|\nabla^{\beta} H\|_{L^2(\mathbb T^3)}^{1-\theta}
\end{equation}
with $q=\frac{6}{7-2\alpha-2\alpha\theta}$ and $p\theta=2$. Taking $q=2+\epsilon$ for some small constant $\epsilon>0$, we obtain $p'=\frac{\alpha}{\alpha-1-\epsilon}$ for a different small constant $\epsilon>0$. It follows that
\begin{equation}\notag
\begin{split}
&\|\nabla^{2-\alpha+\beta} H\|_{L^{p}((0, t_0); L^{q}(\mathbb T^3))}\\
\lesssim &\ \left(\int_0^{t_0} \|\nabla^{\alpha+\beta} H\|_{L^2(\mathbb T^3)}^2 \|\nabla^{\beta}H\|_{L^2(\mathbb T^3)}^{p-2} \, dt\right)^{\frac1p}\\
\lesssim &\ \left(\sup_{0\leq t\leq t_0}\|\nabla^{\beta}H\|_{L^2(\mathbb T^3)}^2 \right)^{\frac{p-2}{2p}}
\left(\int_0^{t_0} \|\nabla^{\alpha+\beta} H\|_{L^2(\mathbb T^3)}^2 \, dt\right)^{\frac{1}{p}}\\
\lesssim&\ \mathcal E_{\alpha}^{\frac12}(H)(t_0).
\end{split}
\end{equation}
Consequently it leads to
\begin{equation}\notag
\|B_{f^\omega} \nabla^{2-\alpha+\beta} H\|_{L^2((0, t_0); L^{2}(\mathbb T^3))}
\lesssim \lambda t_0^{-\gamma}\mathcal E_{\alpha}^{\frac12}(H)(t_0).
\end{equation}
In conclusion we get 
\begin{equation}\label{est-q8}
\|H\otimes B_{f^\omega}\|_{L^2((0, t_0); \mathcal H^{2-\alpha+\beta}(\mathbb T^3))}
\lesssim \lambda t_0^{-\gamma}\mathcal E_{\alpha}^{\frac12}(H)(t_0).
\end{equation}

Thanks to condition (\ref{ass-bf2}), the last term in (\ref{est-q6}) can be estimated
\begin{equation}\label{est-q9}
\begin{split}
&\|B_{f^\omega}\otimes B_{f^\omega}\|_{L^2((0, t_0); \mathcal H^{2-\alpha+\beta}(\mathbb T^3))}\\
\lesssim&\ \|B_{f^\omega} \nabla^{2-\alpha+\beta}B_{f^\omega}\|_{L^2((0, t_0); L^2(\mathbb T^3))}\\
\lesssim&\ \|B_{f^\omega}\|_{L^p((0, t_0); L^{p}(\mathbb T^3))}\|\nabla^{2-\alpha+\beta} B_{f^\omega}\|_{L^{p'}((0, t_0); L^{p'}(\mathbb T^3))}\\
\lesssim&\ \lambda^2t_0^{-2\gamma}
\end{split}
\end{equation}
with $\frac{1}{p}+\frac{1}{p'}=\frac12$.  In order to obtain the largest possible value of $s$ such that for some $\gamma<0$
\[p'\left(\frac{2-\alpha+\beta+s}{2\alpha}-\gamma\right)<1,\ \ \ p\left(\frac{s}{2\alpha}-\gamma\right)<1,\] 
we choose
\[p=\frac{4\alpha}{2\alpha-\beta-2}, \ \ \ p'=\frac{4\alpha}{\beta+2},\]
and hence
\[ s<\frac{2\alpha}{p}+2\alpha\gamma=\frac12(2\alpha-\beta-2)+2\alpha\gamma.\]
Recall that $\beta=\frac72-2\alpha$. Requiring $s>0$ leads to $2\alpha-\beta-2>0$ which implies $\alpha>\frac{11}{8}$.

Therefore, for $\alpha\in(\frac{11}{8}, \frac{7}{4}]$, $0<s<2\alpha-\frac{11}{4}+2\alpha\gamma$ with some $\gamma<0$, we have from (\ref{E2}) and (\ref{est-q6})-(\ref{est-q9}) that
\begin{equation}\notag
\begin{split}
\mathcal E^{\frac12}_{\alpha}(H)(t_0)
\lesssim&\ \mathcal E_{\alpha}(H)(t_0)+\lambda t_0^{-\gamma} \mathcal E^{\frac12}_{\alpha}(H)(t_0)+\lambda^2 t_0^{-2\gamma}, \ \ \mbox{for} \ \ t\in[0,t_0].
\end{split}
\end{equation}
Similarly a continuity argument yields that for small enough $t_0$ we have
\begin{equation}\notag
\mathcal E(H)(t)\leq \mathcal E_{\alpha}(H)(t)\leq C \ \ \ \forall t\in[0,t_0].
\end{equation}

\medskip

{\textbf{(iii) Estimates on $[t_0, T]$ in both 2D and 3D.}}
For $t\in[t_0, T]$, taking inner product of (\ref{H}) with $H$ and integrating over $\mathbb T^n$ yields
\begin{equation}\notag
\begin{split}
&\frac12\frac{d}{dt}\|H(t)\|^2_{L^2(\mathbb T^n)}+\int_{\mathbb T^n} |\nabla^\alpha H|^2\, dx\\
=&-\int_{\mathbb T^n} \left[\nabla\times\nabla\cdot\left((H+B_{f^\omega})\otimes(H+B_{f^\omega})\right)\right]\cdot H\, dx\\
=&-\int_{\mathbb T^n} \left[\nabla\times\nabla\cdot\left((H+B_{f^\omega})\otimes(H+B_{f^\omega})\right)\right]\cdot (H+B_{f^\omega})\, dx\\
&+\int_{\mathbb T^n} \left[\nabla\times\nabla\cdot\left((H+B_{f^\omega})\otimes(H+B_{f^\omega})\right)\right]\cdot B_{f^\omega}\, dx.
\end{split}
\end{equation}
Note that the first integral on the right hand side above is zero due to the fact
\[\int_{\mathbb T^n} \left[\nabla\times\nabla\cdot\left(u\otimes u\right)\right]\cdot u\, dx=\int_{\mathbb T^n} \left[\nabla\times((\nabla\times u)\times u)\right]\cdot u\, dx=0\]
for any vector field $u$ with $\nabla\cdot u=0$. For the same reason, we have 
\[\int_{\mathbb T^n} \left[\nabla\times\nabla\cdot\left(B_{f^\omega}\otimes B_{f^\omega}\right)\right]\cdot B_{f^\omega}\, dx=0.\]
Therefore it follows that
\begin{equation}\label{est-largetime}
\begin{split}
&\frac12\frac{d}{dt}\|H(t)\|^2_{L^2(\mathbb T^n)}+\int_{\mathbb T^n} |\nabla^\alpha H|^2\, dx\\
=&\int_{\mathbb T^n} \left[\nabla\times\nabla\cdot\left((H+B_{f^\omega})\otimes(H+B_{f^\omega})\right)\right]\cdot B_{f^\omega}\, dx\\
=&\int_{\mathbb T^n} \left[\nabla\times\nabla\cdot(H\otimes H)\right]\cdot B_{f^\omega}\, dx\\
&+\int_{\mathbb T^n} \left[\nabla\times\nabla\cdot\left(H\otimes B_{f^\omega}\right)\right]\cdot B_{f^\omega}\, dx\\
&+\int_{\mathbb T^n} \left[\nabla\times\nabla\cdot\left(B_{f^\omega}\otimes H\right)\right]\cdot B_{f^\omega}\, dx\\
=:& \ K_1+K_2+K_3.
\end{split}
\end{equation}
Applying integration by parts we obtain
\begin{equation}\notag
K_1=-\int_{\mathbb T^n} (H\otimes H)\cdot \nabla\nabla\times B_{f^\omega}\, dx.
\end{equation}
It then follows from H\"older's inequality and condition (\ref{ass-bf1})
\begin{equation}\label{est-k1}
\begin{split}
|K_1|\leq&\ \|H\|^2_{L^2(\mathbb T^n)}\|\nabla\nabla\times B_{f^\omega}\|_{L^\infty(\mathbb T^n)}\\
\lesssim&\ (\mathrm{max}\{t^{-\frac12}, t^{-\frac{4+n+2s}{2\alpha}}\})^{\frac12}\|H\|^2_{L^2(\mathbb T^n)}.
\end{split}
\end{equation}
Similarly we have from (\ref{ass-bf1})
\begin{equation}\label{est-k2}
\begin{split}
|K_2|+|K_3|\lesssim&\ \|H\|_{L^2(\mathbb T^n)}\|B_{f^\omega}\|_{L^2}\|\nabla\nabla\times B_{f^\omega}\|_{L^\infty(\mathbb T^n)}\\
\lesssim&\ (1+t^{-\frac{s}{2\alpha}})(\mathrm{max}\{t^{-\frac12}, t^{-\frac{4+n+2s}{2\alpha}}\})^{\frac12}\|H\|_{L^2(\mathbb T^n)}.
\end{split}
\end{equation}
Putting (\ref{est-largetime}), (\ref{est-k1}) and (\ref{est-k2}) together we obtain
\begin{equation}\label{est-large2}
\begin{split}
\frac{d}{dt} \mathcal E(H)(t)\lesssim &\ (\mathrm{max}\{t^{-\frac12}, t^{-\frac{4+n+2s}{2\alpha}}\})^{\frac12}\mathcal E(H)(t)\\
&+(1+t^{-\frac{s}{2\alpha}})(\mathrm{max}\{t^{-\frac12}, t^{-\frac{4+n+2s}{2\alpha}}\})^{\frac12}\mathcal E^{\frac12}(H)(t).
\end{split}
\end{equation}
Note that
\begin{equation}\label{est-time1}
\begin{split}
\int_{t_0}^T(\mathrm{max}\{t^{-\frac12}, t^{-\frac{4+n+2s}{2\alpha}}\})^{\frac12} \, dt=&\int_{t_0}^1t^{-\frac{4+n+2s}{4\alpha}}\,dt+\int_{1}^Tt^{-\frac14}\, dt\\
\leq &\ C(t_0, T, \alpha, n, s)
\end{split}
\end{equation}
and similarly 
\begin{equation}\label{est-time2}
\begin{split}
\int_{t_0}^T(1+t^{-\frac{s}{2\alpha}})(\mathrm{max}\{t^{-\frac12}, t^{-\frac{4+n+2s}{2\alpha}}\})^{\frac12} \, dt
\leq &\ C(t_0, T, \alpha, n, s).
\end{split}
\end{equation}
It follows from (\ref{est-large2}), (\ref{est-time1}) and (\ref{est-time2}) that
\begin{equation}\notag
\mathcal E(H)(t)\leq C(t_0, T, \alpha, n, s) \ \ \forall \ t\in[t_0, T]. 
\end{equation}

\medskip

{\textbf{(iv) Estimates of $\frac{d H}{dt}$ in both 2D and 3D.}} It follows directly from (\ref{H}) that
\begin{equation}\label{est-ht1}
\begin{split}
&\left\|\frac{d}{dt} H\right\|_{L^p([0,T]; \mathcal H^{-2\alpha}(\mathbb T^n))}\\
\lesssim&\ \|(-\Delta)^{\alpha}H\|_{L^p([0,T]; \mathcal H^{-2\alpha}(\mathbb T^n))}+\|\nabla\times \nabla\cdot(H\otimes H)\|_{L^p([0,T]; \mathcal H^{-2\alpha}(\mathbb T^n))}\\
&+\|\nabla\times \nabla\cdot(H\otimes B_{f^\omega})\|_{L^p([0,T]; \mathcal H^{-2\alpha}(\mathbb T^n))}\\
&+\|\nabla\times \nabla\cdot(B_{f^\omega}\otimes B_{f^\omega})\|_{L^p([0,T]; \mathcal H^{-2\alpha}(\mathbb T^n))}.
\end{split}
\end{equation}
When $n=2$, we take $p=2$. It is obvious that
\begin{equation}\notag
\|(-\Delta)^{\alpha}H\|_{L^2([0,T]; \mathcal H^{-2\alpha}(\mathbb T^2))}\lesssim \|H\|_{L^2([0,T]; L^{2}(\mathbb T^2))}\lesssim T^{\frac12}\|H\|_{L^\infty([0,T]; L^{2}(\mathbb T^2))}.
\end{equation}
H\"older's, interpolation and Sobolev embedding inequalities yields for $\alpha>1$
\begin{equation}\notag
\begin{split}
&\|\nabla\times \nabla\cdot(H\otimes H)\|_{L^2([0,T]; \mathcal H^{-2\alpha}(\mathbb T^2))}\\
\lesssim &\ \|H\otimes H\|_{L^2([0, T]; L^{2}(\mathbb T^2))}\\
\lesssim &\ \left(\int_0^{T} \|H\|_{L^2(\mathbb T^2)}^2\|H\|_{L^\infty(\mathbb T^2)}^2\, dt\right)^{\frac12}\\
\lesssim &\ \left(\int_0^{T} \|H\|_{L^2(\mathbb T^2)}^2\|H\|_{\mathcal H^\alpha(\mathbb T^2)}^2\, dt\right)^{\frac12}\\
\lesssim &\ \left(\sup_{t\in(0, T)}\|H(t)\|_{L^2(\mathbb T^2)}\right)\left(\int_0^{T} \|H\|_{\mathcal H^\alpha(\mathbb T^2)}^2\, dt\right)^{\frac12}\\
\lesssim &\ \mathcal E(H)(T).
\end{split}
\end{equation}
It follows from H\"older's inequality and condition (\ref{ass-bf2}) that
for $p$, $p'$ and $m$ satisfying 
\[\frac{1}{p}+\frac{1}{p'}=\frac12, \ \ p'=\frac{4\alpha}{4\alpha-5}, \ \ m=\frac{4\alpha-5}{5-2\alpha}\]
we have
\begin{equation}\notag
\begin{split}
&\|\nabla\times \nabla\cdot(H\otimes B_{f^\omega})\|_{L^2([0,T]; \mathcal H^{-2\alpha}(\mathbb T^2))}\\
\lesssim&\ \|H\otimes B_{f^\omega}\|_{L^2([0, T]; L^{2}(\mathbb T^2))}\\
\lesssim&\ \|H\|_{L^{p}([0,T]; L^{p}(\mathbb T^2))} \|B_{f^\omega}\|_{L^{p'}([0, T]; L^{p'}(\mathbb T^2))}\\
\lesssim&\ \left(\sup_{0\leq t\leq T}\|H(t)\|^2_{L^2(\mathbb T^2)}\right)^{\frac{p-2}{2p}}\|\nabla^mH\|_{L^{2}([0,T]; L^{2}(\mathbb T^2))}^{\frac{2}{p}} \|B_{f^\omega}\|_{L^{p'}([0, T]; L^{p'}(\mathbb T^2))}\\
\lesssim&\ \lambda \mathcal E^{\frac12}(H) T^{-\gamma},
\end{split}
\end{equation}
where we used the fact $m\leq \alpha$ for $\frac43\leq \alpha\leq \frac32$. 
Moreover, the condition (\ref{ass-bf2}) implies
\begin{equation}\notag
\begin{split}
&\|\nabla\times \nabla\cdot(B_{f^\omega}\otimes B_{f^\omega})\|_{L^2([0,T]; \mathcal H^{-2\alpha}(\mathbb T^2))}\\
\lesssim &\ \|B_{f^\omega}\otimes B_{f^\omega}\|_{L^2([0, T]; L^{2}(\mathbb T^2))}\\
\lesssim&\ \|B_{f^\omega}\|^2_{L^4([0, T]; L^4(\mathbb T^2))} \\
\lesssim&\ C(T)\|B_{f^\omega}\|^2_{L^{\frac{4\alpha}{4\alpha-5}}([0, T]; L^{\frac{4\alpha}{4\alpha-5}}(\mathbb T^2))} \\
\leq &\ C(T)\lambda^2T^{-2\gamma}
\end{split}
\end{equation}
since $4\leq \frac{4\alpha}{4\alpha-5}$ for $\frac43\leq \alpha\leq \frac32$.
Combining the estimates above with (\ref{est-ht1}) we have
\begin{equation}\notag
\left\|\frac{d}{dt} H\right\|_{L^2([0,T]; \mathcal H^{-2\alpha}(\mathbb T^2))}\lesssim \mathcal E(H)(T)+\lambda \mathcal E^{\frac12}(H) T^{-\gamma}+\lambda^2T^{-2\gamma}\lesssim C(T, \lambda, s).
\end{equation}

When $n=3$, take $p=\frac{4\alpha}{3}$. First we have
\begin{equation}\notag
\|(-\Delta)^{\alpha}H\|_{L^{\frac{4\alpha}{3}}([0,T]; \mathcal H^{-2\alpha}(\mathbb T^2))}\lesssim \|H\|_{L^{\frac{4\alpha}{3}}([0,T]; L^{2}(\mathbb T^2))}\lesssim T^{\frac{3}{4\alpha}}\|H\|_{L^\infty([0,T]; L^{2}(\mathbb T^2))}.
\end{equation}
Similarly following the application of H\"older's and Gagliardo-Nirenberg's interpolation inequalities we infer
\begin{equation}\notag
\begin{split}
&\|\nabla\times \nabla\cdot(H\otimes H)\|_{L^{\frac{4\alpha}{3}}([0,T]; \mathcal H^{-2\alpha}(\mathbb T^3))}\\
\lesssim &\ \|H\otimes H\|_{L^{\frac{4\alpha}{3}}([0, T];  L^{2}(\mathbb T^3))}\\
\lesssim &\ \left(\int_0^{T} \|\nabla^\alpha H\|_{L^2(\mathbb T^3)}^2\|H\|_{L^2(\mathbb T^3)}^{\frac{8\alpha}{3}-2}\, dt\right)^{\frac{3}{4\alpha}}\\
\lesssim &\ \left(\int_0^{T} \|H\|_{L^2(\mathbb T^3)}^2\|H\|_{\mathcal H^\alpha(\mathbb T^3)}^2\, dt\right)^{\frac{3}{4\alpha}}\\
\lesssim &\ \left(\sup_{t\in(0, T)}\|H(t)\|^2_{L^2(\mathbb T^3)}\right)^{\frac{4\alpha-3}{4\alpha}}\left(\int_0^{T} \|H\|_{\mathcal H^\alpha(\mathbb T^3)}^2\, dt\right)^{\frac{3}{4\alpha}}\\
\lesssim &\ \mathcal E(H)(T).
\end{split}
\end{equation}
For $p, q, p'$ and $m$ satisfying
\[\frac{1}{p}+\frac{1}{p'}=\frac{3}{4\alpha}, \ \ \frac{1}{q}+\frac{1}{p'}=\frac{1}{2}, \ \ p'=\frac{8\alpha}{8\alpha-11}, \ \ m=\frac{3(8\alpha-11)}{2(17-8\alpha)},\]
we apply H\"older's inequality and condition (\ref{ass-bf2}) to deduce
\begin{equation}\notag
\begin{split}
&\|\nabla\times \nabla\cdot(H\otimes B_{f^\omega})\|_{L^{\frac{4\alpha}{3}}([0,T]; \mathcal H^{-2\alpha}(\mathbb T^3))}\\
\lesssim&\ \|H\otimes B_{f^\omega}\|_{L^{\frac{4\alpha}{3}}([0, T]; L^{2}(\mathbb T^3))}\\
\lesssim&\ \|H\|_{L^{p}([0,T]; L^{q}(\mathbb T^3))} \|B_{f^\omega}\|_{L^{p'}([0, T]; L^{p'}(\mathbb T^3))}\\
\lesssim&\ \left(\sup_{0\leq t\leq T}\|H(t)\|^2_{L^2(\mathbb T^3)}\right)^{\frac{p-2}{2p}}\|\nabla^mH\|_{L^{2}([0,T]; L^{2}(\mathbb T^3))}^{\frac{2}{p}} \|B_{f^\omega}\|_{L^{p'}([0, T]; L^{p'}(\mathbb T^3))}\\
\lesssim&\ \lambda \mathcal E^{\frac12}(H) T^{-\gamma}
\end{split}
\end{equation}
thanks to the fact that $m\leq \alpha$ for $\frac{11}{8}<\alpha\leq \frac{7}{4}$.
In the end, it follows from H\"older's inequality and (\ref{ass-bf2}) that
\begin{equation}\notag
\begin{split}
&\|\nabla\times \nabla\cdot(B_{f^\omega}\otimes B_{f^\omega})\|_{L^{\frac{4\alpha}{3}}([0,T]; \mathcal H^{-2\alpha}(\mathbb T^3))}\\
\lesssim &\ \|B_{f^\omega}\otimes B_{f^\omega}\|_{L^{\frac{4\alpha}{3}}([0, T]; L^{2}(\mathbb T^3))}\\
\lesssim&\ \|B_{f^\omega}\|^2_{L^{\frac{8\alpha}{3}}([0, T]; L^4(\mathbb T^3))} \\
\lesssim&\ C(T)\|B_{f^\omega}\|^2_{L^{\frac{8\alpha}{8\alpha-11}}([0, T]; L^{\frac{8\alpha}{8\alpha-11}}(\mathbb T^2))} \\
\leq &\ C(T)\lambda^2T^{-2\gamma}
\end{split}
\end{equation}
since $\frac{8\alpha}{3}\leq \frac{8\alpha}{8\alpha-11}$ for $\alpha\leq \frac74$.
Again collecting the estimates above with (\ref{est-ht1}) we get
\begin{equation}\notag
\left\|\frac{d}{dt} H\right\|_{L^{\frac{4\alpha}{3}}([0,T]; \mathcal H^{-2\alpha}(\mathbb T^3))}
\lesssim C(T, \lambda, s).
\end{equation}

\cbdu

\bigskip

\section{Existence of weak solutions to (\ref{H})}\label{sec-galerkin}

We are ready to establish the existence of weak solutions to (\ref{H}) by using the standard Galerkin approximating approach (c.f. \cite{CF, DG}) and the a priori estimates obtained in the previous section. Namely we will prove Theorem \ref{thm-H} by constructing a sequence of Galerkin approximating solutions and passing to a limit. 

Recall the Fourier transform and its inverse on torus $\mathbb T^n$,
\begin{equation}\notag
\begin{split}
\widehat f(k, t)=&\int_{\mathbb T^n} f(x,t)e^{-2\pi ik\cdot x}\, dx, \ \ \ \ k\in \mathbb Z^n\\
f(x,t)=&\sum_{k\in\mathbb Z^n} \widehat f(k, t)e^{2\pi ik\cdot x}.
\end{split}
\end{equation}
Denote $P_K$ by the Fourier projection operator
\begin{equation}\notag
P_K f=\sum_{\{k: |k_i|\leq K, 1\leq i\leq n\}} \widehat f(k, t)e^{2\pi ik\cdot x}
\end{equation}
and $H^K= P_K H$. 
For any fixed $K\in \mathbb N$ we consider the truncated system 
\begin{equation}\label{eq-K}
\begin{split}
H^K_t=&-(-\Delta)^\alpha H^K-P_K\left[\nabla\times\left(\mathcal B(H^K, H^K)+\mathcal B(H^K, B^K_{f^\omega}) \right)\right]\\
&-P_K\left[\nabla\times \left(\mathcal B(B^K_{f^\omega}, H^K)+\mathcal B(B^K_{f^\omega}, B^K_{f^\omega})\right)\right],\\
\nabla\cdot H^K=&\ 0,\\
H^K(x,0)=&\ 0. 
\end{split}
\end{equation}
Taking Fourier transform on (\ref{eq-K}) yields
\begin{equation}\label{eq-KF}
\begin{split}
\widehat {H^K}_t=&\ (-1)^\alpha |k|^{2\alpha}\widehat {H^K}(k,t)\\
&-ik\times \sum_{\{k'+k''=k, |k'_i|\leq K, |k''_i|\leq K \}}\widehat {H^K} (k', t)\cdot k'' \widehat{H^K}(k'', t)\\
&-ik\times \sum_{\{k'+k''=k, |k'_i|\leq K, |k''_i|\leq K \}}\widehat {H^K} (k', t)\cdot k'' \widehat{B^K_{f^\omega}}(k'', t)\\
&-ik\times \sum_{\{k'+k''=k, |k'_i|\leq K, |k''_i|\leq K \}}\widehat {B^K_{f^\omega}} (k', t)\cdot k'' \widehat{H^K}(k'', t)\\
&-ik\times \sum_{\{k'+k''=k, |k'_i|\leq K, |k''_i|\leq K \}}\widehat {B^K_{f^\omega}} (k', t)\cdot k'' \widehat{B^K_{f^\omega}}(k'', t),\\
k\cdot \widehat{H^K}(k, t)=&\ 0,\\
\widehat{H^K}(k,0)=&\ 0. 
\end{split}
\end{equation}
Note that (\ref{eq-KF}) is a finite ODE system for any fixed $K\in\mathbb N$. From the integral form of (\ref{eq-KF}) we define the map
\begin{equation}\notag
\begin{split}
\Phi(\widehat{H^K})(k,t):=& \int_0^t(-1)^\alpha |k|^{2\alpha}\widehat {H^K}(k,s)\, ds\\
&-\int_0^t ik\times \sum_{\{k'+k''=k, |k'_i|\leq K, |k''_i|\leq K \}}\widehat {H^K} (k', t)\cdot k'' \widehat{H^K}(k'', s)\, ds\\
&-\int_0^t ik\times \sum_{\{k'+k''=k, |k'_i|\leq K, |k''_i|\leq K \}}\widehat {H^K} (k', t)\cdot k'' \widehat{B^K_{f^\omega}}(k'', s)\, ds\\
&-\int_0^t ik\times \sum_{\{k'+k''=k, |k'_i|\leq K, |k''_i|\leq K \}}\widehat {B^K_{f^\omega}} (k', t)\cdot k'' \widehat{H^K}(k'', s)\, ds\\
&-\int_0^t ik\times \sum_{\{k'+k''=k, |k'_i|\leq K, |k''_i|\leq K \}}\widehat {B^K_{f^\omega}} (k', t)\cdot k'' \widehat{B^K_{f^\omega}}(k'', s)\, ds\\
=:&\ \Phi_1(k,t)+\Phi_2(k,t)+\Phi_3(k,t)+\Phi_4(k,t)+\Phi_5(k,t).
\end{split}
\end{equation}
Denote the function space
\begin{equation}\notag
X_T= C([0,T]; \ell^2)\cap L^2([0,T]; \mathcal H^\alpha), \ \ \mbox{for} \ \ T>0.
\end{equation}
We first show that the map $\Phi$ has a fixed point on $X_{t_1}$ for a small time $t_1$ by showing that $\Phi$ is a contraction map on a ball of $X_{t_1}$. We then claim that this process can be iterated to reach time $T$. 

For $t\in[0,t_1]$ one has
\begin{equation}\notag
\|\Phi_1(t)\|_{\ell^2}\lesssim K^{2\alpha} t_1\|\widehat{H^K}\|_{L^\infty([0,t_1]; \ell^2)}.
\end{equation} 
Applying Plancherel's theorem and Sobolev imbedding gives
\begin{equation}\notag
\|\Phi_2(t)\|_{\ell^2}\lesssim K^{2+\frac{n}{2}} t_1 \|\widehat{H^K}\|^2_{L^\infty([0,t_1]; \ell^2)}.
\end{equation} 
Using Plancherel's theorem and Sobolev imbedding again and the estimate (\ref{est-lin1}) we have
\begin{equation}\notag
\|\Phi_3(t)\|_{\ell^2}+\|\Phi_4(t)\|_{\ell^2} \lesssim K^{2+\frac{n}{2}} t_1^{1-\frac{s}{2\alpha}}\|\widehat{H^K}\|_{L^\infty([0,t_1]; \ell^2)}.
\end{equation} 
 It follows from Plancherel's theorem and Lemmas \ref{le-heat3} and \ref{le-heat4} that
 \begin{equation}\notag
\|\Phi_5(t)\|_{\ell^2}\lesssim K^2\lambda^2 t_1^{\frac{\beta+2}{2\alpha}-2\gamma},
\end{equation} 
where we recall $\beta=3-2\alpha$ in 2D and $\beta=\frac72-2\alpha$ in 3D, and we observe $\frac{\beta+2}{2\alpha}>0$.
Thus combining the estimates above leads to
\begin{equation}\label{est-contract1}
\begin{split}
\|\Phi(\widehat{H^K})(t)\|_{\ell^2}\lesssim&\ K^{2\alpha} t_1\|\widehat{H^K}\|_{L^\infty([0,t_1]; \ell^2)}+K^{2+\frac{n}{2}} t_1 \|\widehat{H^K}\|^2_{L^\infty([0,t_1]; \ell^2)}\\
&+K^{2+\frac{n}{2}} t_1^{1-\frac{s}{2\alpha}}\|\widehat{H^K}\|_{L^\infty([0,t_1]; \ell^2)}+K^2\lambda^2 t_1^{\frac{\beta+2}{2\alpha}-2\gamma}. 
\end{split}
\end{equation} 
Analogously we obtain
\begin{equation}\label{est-contract2}
\begin{split}
&\||k|^\alpha \Phi(\widehat{H^K})(t)\|_{L^2([0,t_1];\ell^2)}\\
\lesssim&\ K^{3\alpha} t_1\|\widehat{H^K}\|_{L^\infty([0,t_1]; \ell^2)}+K^{2+\alpha+\frac{n}{2}} t_1 \|\widehat{H^K}\|^2_{L^\infty([0,t_1]; \ell^2)}\\
&+K^{2+\alpha+\frac{n}{2}} \lambda t_1^{1-\frac{s}{2\alpha}}\|\widehat{H^K}\|_{L^\infty([0,t_1]; \ell^2)}+K^{2+\alpha}\lambda^2 t_1^{\frac{\beta+2}{2\alpha}-2\gamma}. 
\end{split}
\end{equation} 
Take $R=K$. Note that $\gamma<0$, $s<2\alpha$ and $\frac{\beta+2}{2\alpha}>0$; hence all the index of $t_1$ are positive in (\ref{est-contract1})-(\ref{est-contract2}).  Thus we can choose $t_1$ small enough such that the estimates (\ref{est-contract1}) and (\ref{est-contract2}) imply that $\Phi$ maps the ball $B_R(0)\subset X_{t_1}$ to itself continuously. Analogous analysis guarantees that the map $\Phi$ is a contraction. Hence there exists a unique solution $\widehat{H^K}$ to (\ref{eq-KF}) in $X_{t_1}$. Consequently, there exists a unique solution $H^K$ to (\ref{eq-K}) in $L^\infty([0,t_1]; L^2(\mathbb T^n))\cap L^2([0,t_1]; \mathcal H^\alpha(\mathbb T^n))$. Note that since the energy estimate (\ref{ap-est1}) holds for system (\ref{eq-K}) on $[0,T]$ as well, iterations of the previous process can yield a solution of (\ref{eq-K}) up to time $T$. Automatically the solution $H^K$ satisfies estimates (\ref{ap-est1}) and (\ref{ap-est2}). Note that $P_K$ is a bounded operator in $L^p$ for any $1<p<\infty$ and hence $B^K_{f^\omega}$ converges strongly to $B_{f^\omega}$ in $L^p$ as $K\to\infty$. Therefore the estimates (\ref{ap-est1}) and (\ref{ap-est2}) are sufficient for us to extract a subsequence of $H^K$ which converges to a weak solution $H$ of (\ref{H}) on $[0,T]$.

\bigskip

\section{Appendix: Proof of uniqueness}\label{sec-unique}
In this section, we show the uniqueness of the weak solutions for critical and subcritical values of $\alpha$. Let $B^1=B_{f^\omega}+H^1$ and  $B^2=B_{f^\omega}+H^2$ be two weak solutions obtained in Theorems \ref{thm-2d} and \ref{thm-3d} in 2D and 3D respectively, for system (\ref{emhd}) with the same initial data $f$. Thus both $H^1$ and $H^2$ satisfy (\ref{H}). In order to fully explore cancellations in the estimates later, we write the equations of $H^1$ and $H^2$ as 
\begin{equation}\notag
\begin{split}
H^1_t+\nabla\times((\nabla\times(H^1+B_{f^\omega}))\times (H^1+B_{f^\omega}))=-(-\Delta)^\alpha H^1,\\
H^2_t+\nabla\times((\nabla\times(H^2+B_{f^\omega}))\times (H^2+B_{f^\omega}))=-(-\Delta)^\alpha H^2.
\end{split}
\end{equation}
Denote $\widetilde H=H^1-H^2$. Taking subtraction of the last two equations gives
\begin{equation}\label{eq-diff}
\widetilde H_t+\nabla\times((\nabla\times\widetilde H)\times (H^1+B_{f^\omega}))+\nabla\times((\nabla\times(H^2+B_{f^\omega}))\times \widetilde H)=-(-\Delta)^\alpha \widetilde H.
\end{equation}
Taking inner product of (\ref{eq-diff}) with $\widetilde H$, integrating over $\mathbb T^n$ and using integration by parts we obtain
\begin{equation}\label{energy-diff1}
\begin{split}
&\frac{1}{2}\frac{d}{dt} \|\widetilde H(t)\|_{L^2(\mathbb T^n)}^2+\int_{\mathbb T^n}|\nabla^\alpha \widetilde H(t)|^2\, dx\\
=& -\int_{\mathbb T^n} (\nabla\times(H^2(t)+B_{f^\omega}))\times \widetilde H(t)\cdot \nabla\times \widetilde H(t)\, dx
\end{split}
\end{equation}
where we used the cancellation
\begin{equation}\notag
\begin{split}
\int_{\mathbb T^n}\nabla\times((\nabla\times\widetilde H)\times (H^1+B_{f^\omega}))\cdot  \widetilde H\, dx
=&\int_{\mathbb T^n}((\nabla\times\widetilde H)\times (H^1+B_{f^\omega}))\cdot \nabla\times \widetilde H\, dx\\
=&\ 0.
\end{split}
\end{equation}
In 2D, i.e. $n=2$, we estimate the integral on the right hand side of (\ref{energy-diff1}) by using H\"older's inequality, Sobolev embedding and Young's inequality for some $p$ and $q$ satisfying $\frac{1}{p}+\frac{1}{q}=\frac12$
\begin{equation}\notag
\begin{split}
&\left| \int_{\mathbb T^2} (\nabla\times H^2(t))\times \widetilde H(t)\cdot \nabla\times \widetilde H(t)\, dx\right|\\
\leq &\ C\|\widetilde H\|_{L^2(\mathbb T^2)} \|\nabla H^2\|_{L^p(\mathbb T^2)}  \|\nabla \widetilde H\|_{L^q(\mathbb T^2)} \\
\leq &\ C\|\widetilde H\|_{L^2(\mathbb T^2)} \|\nabla^{2-\frac{2}{p}} H^2\|_{L^2(\mathbb T^2)}  \|\nabla^{2-\frac{2}{q}} \widetilde H\|_{L^2(\mathbb T^2)} \\
\leq &\ C\|\widetilde H\|_{L^2(\mathbb T^2)} \|\nabla^{\alpha} H^2\|_{L^2(\mathbb T^2)}  \|\nabla^{\alpha} \widetilde H\|_{L^2(\mathbb T^2)} \\
\leq &\ C\|\widetilde H\|^2_{L^2(\mathbb T^2)} \|\nabla^{\alpha} H^2\|^2_{L^2(\mathbb T^2)} +\frac12 \|\nabla^{\alpha} \widetilde H\|^2_{L^2(\mathbb T^2)} \\
\end{split}
\end{equation}
where we require $2-\frac{2}{p}\leq \alpha$ and $2-\frac{2}{q}\leq \alpha$. When $\alpha\geq \frac32$, we are able to find proper $p$ and $q$ satisfying the conditions. Analogously, in 3D we have for $\alpha\geq \frac74$
\begin{equation}\notag
\begin{split}
&\left| \int_{\mathbb T^3} (\nabla\times H^2(t))\times \widetilde H(t)\cdot \nabla\times \widetilde H(t)\, dx\right|\\
\leq &\ C\|\widetilde H\|_{L^2(\mathbb T^3)} \|\nabla H^2\|_{L^p(\mathbb T^3)}  \|\nabla \widetilde H\|_{L^q(\mathbb T^3)} \\
\leq &\ C\|\widetilde H\|_{L^2(\mathbb T^3)} \|\nabla^{\frac52-\frac{3}{p}} H^2\|_{L^2(\mathbb T^3)}  \|\nabla^{\frac52-\frac{3}{q}} \widetilde H\|_{L^2(\mathbb T^3)} \\
\leq &\ C\|\widetilde H\|^2_{L^2(\mathbb T^3)} \|\nabla^{\alpha} H^2\|^2_{L^2(\mathbb T^3)} +\frac12 \|\nabla^{\alpha} \widetilde H\|^2_{L^2(\mathbb T^3)} \\
\end{split}
\end{equation}
provided $\frac{1}{p}+\frac{1}{q}=\frac12$, $\frac52-\frac{3}{p}\leq \alpha$ and $\frac52-\frac{3}{q}\leq \alpha$. 

On the other hand, we have in 2D
\begin{equation}\notag
\begin{split}
&\left| \int_{\mathbb T^2} (\nabla\times B_{f^\omega})\times \widetilde H(t)\cdot \nabla\times \widetilde H(t)\, dx\right|\\
\leq &\ C\|\widetilde H\|_{L^2(\mathbb T^2)} \|\nabla B_{f^\omega}\|_{L^p(\mathbb T^2)}  \|\nabla \widetilde H\|_{L^q(\mathbb T^2)} \\
\leq &\ C\|\widetilde H\|_{L^2(\mathbb T^2)} \|\nabla^{2-\frac{2}{p}} B_{f^\omega}\|_{L^2(\mathbb T^2)}  \|\nabla^{2-\frac{2}{q}} \widetilde H\|_{L^2(\mathbb T^2)} \\
\leq &\ C\|\widetilde H\|_{L^2(\mathbb T^2)} \|\nabla^{\alpha} B_{f^\omega}\|_{L^2(\mathbb T^2)}  \|\nabla^{\alpha} \widetilde H\|_{L^2(\mathbb T^2)} \\
\leq &\ C\|\widetilde H\|^2_{L^2(\mathbb T^2)} \|\nabla^{\alpha} B_{f^\omega}\|^2_{L^2(\mathbb T^2)} +\frac12 \|\nabla^{\alpha} \widetilde H\|^2_{L^2(\mathbb T^2)}
\end{split}
\end{equation}
and in 3D
\begin{equation}\notag
\begin{split}
&\left| \int_{\mathbb T^3} (\nabla\times  B_{f^\omega}(t))\times \widetilde H(t)\cdot \nabla\times \widetilde H(t)\, dx\right|\\
\leq &\ C\|\widetilde H\|_{L^2(\mathbb T^3)} \|\nabla  B_{f^\omega}\|_{L^p(\mathbb T^3)}  \|\nabla \widetilde H\|_{L^q(\mathbb T^3)} \\
\leq &\ C\|\widetilde H\|_{L^2(\mathbb T^3)} \|\nabla^{\frac52-\frac{3}{p}}  B_{f^\omega}\|_{L^2(\mathbb T^3)}  \|\nabla^{\frac52-\frac{3}{q}} \widetilde H\|_{L^2(\mathbb T^3)} \\
\leq &\ C\|\widetilde H\|^2_{L^2(\mathbb T^3)} \|\nabla^{\alpha} B_{f^\omega}\|^2_{L^2(\mathbb T^3)} +\frac12 \|\nabla^{\alpha} \widetilde H\|^2_{L^2(\mathbb T^3)} \\
\end{split}
\end{equation}

Therefore, it follows from (\ref{energy-diff1}) that
\begin{equation}\label{energy-diff2}
\begin{split}
&\frac{d}{dt} \|\widetilde H(t)\|_{L^2(\mathbb T^n)}^2+\int_{\mathbb T^n}|\nabla^\alpha \widetilde H(t)|^2\, dx\\
\leq &\ C\|\widetilde H\|^2_{L^2(\mathbb T^n)} \left(\|\nabla^{\alpha} H^2\|^2_{L^2(\mathbb T^n)}+\|\nabla^{\alpha} B_{f^\omega}\|^2_{L^2(\mathbb T^n)} \right).
\end{split}
\end{equation}
Applying Gr\"onwall's inequality to (\ref{energy-diff2}) we obtain
\begin{equation}\label{energy-diff3}
\begin{split}
&\|\widetilde H(t)\|_{L^2(\mathbb T^n)}^2\\
\leq &\ \|\widetilde H(0)\|_{L^2(\mathbb T^n)}^2 \exp\left\{C\int_0^t\left(\|\nabla^{\alpha} H^2\|^2_{L^2(\mathbb T^n)}+\|\nabla^{\alpha} B_{f^\omega}\|^2_{L^2(\mathbb T^n)} \right) \, d\tau \right\}.
\end{split}
\end{equation}
Note that $H^2\in L^2([0,T]; \mathcal H^\alpha (\mathbb T^n))$ and from (\ref{est-lin1})
\begin{equation}\notag
 \int_0^t\|\nabla^{\alpha} B_{f^\omega}\|^2_{L^2(\mathbb T^n)} \, d\tau\lesssim  \int_0^t  \tau^{-\frac{\alpha+s}{2\alpha}}\, d\tau\lesssim t^{1-\frac{\alpha+s}{2\alpha}}.
\end{equation}
Recall that $\alpha\in [\frac43, \frac32]$ and $s\in(0, 2\alpha-\frac52)$ in 2D, and $\alpha\in (\frac{11}{8}, \frac74]$ and $s\in (0, 2\alpha-\frac{11}{4})$ in 3D. One can check that $\frac12<\frac{\alpha+s}{2\alpha}<1$ in both cases. 
Combining with the fact $\widetilde H(0)=0$, (\ref{energy-diff3}) implies
$\widetilde H(t)\equiv 0$ and hence $H^1(t)\equiv H^2(t)$. It follows naturally $B^1(t)\equiv B^2(t)$.

\bigskip

\section*{Acknowledgement}
The author is partially supported by the NSF grants DMS-1815069 and DMS-2009422, and the von Neumann Fellowship.  She is also grateful to IAS for its hospitality in 2021-2022.

\bigskip



\begin{thebibliography}{XX}



\bibitem{ADFL}
M. Acheritogaray, P. Degond, A. Frouvelle and J-G. Liu.
\newblock {\em Kinetic formulation and global existence for the Hall-Magnetohydrodynamic system}.
\newblock Kinetic and Related Models, 4: 901--918, 2011.







\bibitem{Bou96}
J. Bourgain.
\newblock {\em Invariant measures for the 2D defocusing nonlinear Schr\"odinger equation}.
\newblock Comm. Math. Phys., 176: 421-445, 1996.









\bibitem{BT1}
N. Burq and N. Tzvetkov.
\newblock {\em Random data Cauchy theory for super-critical wave equation I: Local theory}.
\newblock Invent. Math., 173(3):449-475, 2008.

\bibitem{BT2}
N. Burq and N. Tzvetkov.
\newblock {\em Random data Cauchy theory for super-critical wave equation II: A global existence result}.
\newblock Invent. Math., 173(3):477-496, 2008.

\bibitem{CDL}
D. Chae, P. Degond and J-G. Liu.
\newblock {\em Well-posedness for Hall-magnetohydrodynamics}.
\newblock Ann. Inst. H. Poincar\'e Anal. Non Lineaire, Vol. 31: 555--565, 2014.



\bibitem{CWW}
D. Chae,  R. Wan and J. Wu.
\newblock {\em Local well-posedness for the Hall--MHD equations with fractional magnetic diffusion}.
\newblock J. Math. Fluid Mech., 17: 627-638, 2015.

\bibitem{CWeng}
D. Chae and S. Weng.
\newblock {\em Singularity formation for the incompressible Hall-MHD equations without resistivity}.
\newblock Ann. I. H. Poincar\'e-AN, Vol. 33: 1009--1022, 2016.

\bibitem{CW}
D. Chae and J. Wolf.
\newblock {\em On partial regularity for the 3D non-stationary Hall magnetohydrodynamics equations on the plane}.
\newblock Comm. Math. Phys., Vol. 354: 213--230, 2017.















\bibitem{CF}
P. Constantin and C. Foias.
\newblock {\em Naiver-Stokes Equations}.
\newblock Chicago Lectures in Mathematics, The University of Chicago Press, 1988.




\bibitem{Dai2}
M. Dai.
\newblock {\em Local well-posedness for the Hall-MHD system in optimal Sobolev spaces}.
\newblock Journal of Differential Equations, 289: 159-181, 2021.














\bibitem{Deng}
Y. Deng.
\newblock {\em Two-dimensional nonlinear Schr\"ondinger equation with random radial data}.
\newblock Anal. PDE, 5: 913-960, 2012.




\bibitem{DG}
C. Doering and J. D. Gibbon.
\newblock {\em Applied Analysis of the Navier-Stokes Equations}.
\newblock Cambridge Texts in Applied Mathematics, Cambridge University Press, 2004.























\bibitem{LR}
P. G. Lemari\'e-Rieusset.
\newblock {\em Recent developments in the Navier-{S}tokes problem}.
\newblock Chapman and Hall/CRC Research Notes in Mathematics, 431. Chapman  and Hall/CRC, Boca Raton, FL, 2002.



\bibitem{LM}
J. L\"uhrmann and D. Mendelson.
\newblock {\em Random data Cauchy theory for nonlinear wave equations of power-type on $\mathbb R^3$}.
\newblock Comm. Partial Differential Equations, 39(12): 2262-2283, 2014.








\bibitem{NPS}
A. R. Nahmod, N. Pavlovi\'c, and G. Staffilani.
\newblock {\em Almost sure existence of global weak solutions for super-critical Navier-Stokes equations}.
\newblock SIAM J. Math. Anal., 45(6): 3431-3452, 2013.












\bibitem{Tao}
T. Tao.
\newblock {\em A quantitative formulation of the global regularity problem for the periodic Navier-Stokes equation}.
\newblock Dyn. Partial Differ. Eq., 4(4): 293-302, 2007.






\end{thebibliography}
\end{document}